\documentclass[12pt,a4paper]{article}
\usepackage{amsfonts, amsopn, amsmath,amssymb,amsthm}
\usepackage[latin1]{inputenc}

\input epsf

\def\supp{\mathop {\rm Supp}}

\def\ad{\mathop {\rm ad}}

\def\Z{\ensuremath{\mathbb Z}}
\def\C{\ensuremath{\mathbb C}}
\def\P{\ensuremath{\mathbb P}}

\newtheorem{theorem}{Theorem}[section]
\newtheorem{corollary}[theorem]{Corollary}
\newtheorem{proposition}[theorem]{Proposition}
\newtheorem{lemma}[theorem]{Lemma}

\theoremstyle{definition}
\newtheorem{definition}{Definition}[section]

\newtheorem*{aknow}{Acknowledgements}
\newtheorem*{notation}{Notation}

\theoremstyle{remark}
\newtheorem{remark}{Remark}[section]
\newtheorem*{howtoreadtableone}{How to read Table\texttilde\ref{tab:weights}}

\numberwithin{equation}{section}

\newcommand{\texttilde}{~}
\newcommand{\Hilb}{\mathrm{Hilb}}

\newcounter{primitivenumber}
\newcommand{\primitivecase}{\smallbreak\noindent\addtocounter{primitivenumber}{1}(\theprimitivenumber)\enspace}

\title{Equivariant deformations of the affine multicone over a flag variety}

\author{P.~Bravi and S.~Cupit-Foutou}

\date{}

\begin{document}

\maketitle

\begin{abstract}
We prove that the invariant Hilbert scheme parametrising the equivariant
deformations of the affine multicone over a flag variety is,
under certain hypotheses, an affine space. 
The proof is based on the construction of a
wonderful variety in a fixed multiprojective space.
\end{abstract}


\section*{Introduction}
Given a finite dimensional $G$-module $V$ with $G$ a complex connected
reductive group,
the invariant Hilbert scheme parameterises the
closed $G$-subschemes of $V$ whose coordinate ring regarded as a $G$-module is isomorphic to a fixed $G$-module.
This fixed $G$-module defines an invariant -- the Hilbert function--  which is the analog of the Hilbert polynomial in case of Grothendieck's Hilbert scheme.
In \cite{AB}, V. Alexeev and M. Brion studied the invariant Hilbert scheme and proved that it was quasiprojective.

Of special interest (\cite{AB}) are the invariant Hilbert schemes whose 
corresponding Hilbert function takes $0$ or $1$ as only values, namely the
fixed $G$-module is multiplicity-free.
Their closed points are the so called affine spherical varieties
(e.g.\ affine toric varieties).

The coordinate ring of the affine cone of primitive vectors in a simple
$G$-module $V$ is multiplicity-free. The equivariant deformations of such a
cone are parameterised by the invariant Hilbert scheme for the corresponding 
Hilbert function.
In \cite{J}, S.~Jansou proved that this invariant Hilbert scheme
was a (reduced) point or an affine line.

In this paper we consider the more general situation 
given by an affine multicone over a flag variety,
provided it is normal, spherical and with boundary of codimension at least
two.

These hypotheses we make on the Hilbert function are not new in invariant
deformation theory (\emph{see}~\cite{Pa} and also \cite{VP,BK,Bri86,K90}).
They yield a correspondence between deformations
of spherical $G$-varieties and deformations of their open $G$-orbits.

Under these assumptions, we thus prove that 
the invariant Hilbert scheme, parameterising the equivariant
deformations of an affine multicone, is an affine space. 
Its dimension can be explicitly determined: it is closely related to a
combinatorial datum (called spherical system) that we can attach to the
Hilbert function.

In \cite{J}, the cases where the Hilbert scheme is the affine line are 
explicitly listed. This list corresponds to that of the two-orbit varieties 
whose closed orbit is a very ample divisor 
(\cite{A} and \cite{HS,Bri89}).

The two-orbit varieties which appear in the above picture are wonderful 
varieties.
The wonderful varieties are spherical (\cite{L96}), the knowledge of
this special class of varieties should lead to a better understanding of 
spherical varieties (\emph{see}~\cite{L01}).
In general, to any wonderful variety one can associate a spherical system,
but the converse remains an open problem; only partial
results are known (loc.\ cit.\ and \cite{Pe03,BP,Bra}).

This paper is an attempt to clarify the connection between invariant
Hilbert schemes and wonderful varieties.
 
In this setting we prove that the total cone of deformations of the
multicone over the flag variety we consider is given by the multicone over
a well determined wonderful $G$-variety containing that flag variety as its
unique closed $G$-orbit.


In the first section we recall the definition of
the invariant Hilbert scheme in the multiplicity free case.
We focus on its closed points and define the
action of the adjoint torus.
We also briefly recall some definitions concerning wonderful varieties.

In the second section we introduce our affine multicone over a flag variety
under certain assumptions and we state our main results.

The remaining sections are dedicated to the proofs.

By means of a  representation theoretical characterisation of the tangent space
of the invariant Hilbert scheme at the multicone
given in \cite{AB}, 
we compute explicitly in the third section the weights of the toric action
on this tangent space.  
We prove also that these weights are multiplicity free.



In the fourth section we recall the definition of spherical systems, due to
D.~Luna. We then
associate a spherical system to the considered invariant Hilbert function.

Afterward, we show that a wonderful spherical subgroup $N$
corresponding to such a spherical system provides an open orbit in the
Hilbert scheme. Furthermore, by means of a result of G.~Pezzini, we prove that
the compactification of $G/N$ in the given
multiprojective space is the investigated wonderful variety.

In the last section we show how to associate to a spherical system in our
setting a corresponding wonderful subgroup.

\begin{aknow}
This work is mostly the result of a staying of the first author at the
Bergische Universit\"at Wuppertal while supported by an European ``Marie
 Curie'' TMR \emph{Liegrits} fellowship.
Both authors would like to thank M.~Brion, S.~Jansou, P.~Littelmann and
D.~Luna for helpful discussions, suggestions and support. 
\end{aknow}

\begin{notation}
The ground field is the field of complex numbers.
Given a connected, simply connected, semisimple algebraic group $G$ and a
maximal torus $T$, 
let $\Lambda$ be the lattice of $T$-weights and $\Phi$ its root system.
Fixing a Borel subgroup $B$ containing $T$, we denote by
$S=\{\alpha_1,\ldots,\alpha_n\}$ the associated set of simple roots 
and by $\Phi^+$ the set of positive roots.
To each root $\beta$, we attach a root vector $X_\beta$ of the Lie algebra
$\mathfrak g$ of $G$. 

Let $\Lambda^+$ be the cone of dominant weights.
Recall that $\Lambda^+=\oplus_{1}^{n}\mathbb N\omega_i$
where $\omega_i$ is the fundamental weight attached to $\alpha_i$.
If $\alpha^\vee _j$ is the coroot of $\alpha_j$,
we have $(\omega_i,\alpha^\vee_j)$ equal to $1$ if $i=j$ and equal to $0$
otherwise. 

The simple $G$-modules are in correspondence with the cone $\Lambda^+$ of
dominant weights. 
We thus denote by $V(\lambda)$ the simple $G$-module with highest weight
$\lambda\in\Lambda^+$. 
\end{notation}

\section{Definitions}

\subsection{Invariant Hilbert scheme}

Let $\Gamma$ be a submonoid of $\Lambda^+$ generated by dominant weights
$\lambda_1,\ldots,\lambda_s$. 
Denote by $V$ the finite dimensional $G$-module given by
$\oplus_{i=1}^sV(\lambda_i)$, 
and by $R$ the $G$-module given by $\oplus_{\lambda\in\Gamma}V(\lambda)^\ast$.

In this setting, we follow \cite{AB} to define (and state some properties
of) the \textit{invariant Hilbert scheme} $\mathrm{Hilb}_\Gamma^G(V)$ of
affine $G$-subschemes of $V$ with $G$-module of regular functions
isomorphic to $R$. In a more general context, the invariant Hilbert scheme
is defined 
up to the \emph{Hilbert function} instead of the weight monoid $\Gamma$. In
the multiplicity free case, these two data are equivalent so we will not
recall here the definition of this function. 

A \textit{family} of affine $G$-subschemes of $V$ over a given scheme $S$
is a closed $G$-subscheme $\mathfrak X\subset V\times S$ where the
restriction on $\mathfrak X$ of the projection $\pi\colon V\times S\to S$
is affine and of finite type. 

Let $\mathcal R$ denote the direct image $\pi_\ast\mathcal O_\mathfrak X$
of the structure sheaf of $\mathfrak X$. Such a family is said to be
\textit{of type $\Gamma$} if we have 
\begin{equation}
\mathcal R\cong
\bigoplus_{\lambda\in\Gamma}\mathcal R_{\lambda}\otimes V(\lambda)^\ast
\end{equation}
and each $\mathcal R_{\lambda}$ is an invertible sheaf of $\mathcal
O_S$-modules. In particular the map $\pi\colon\mathfrak X\to S$ is flat. 

\begin{theorem}[\cite{AB} Theorem~1.7]
The contravariant functor $Hilb_\Gamma^G(V)\colon$ $\{schemes\}\to
\{sets\}$, that to any scheme $S$ associates the set of families $\mathfrak
X\subset V\times S$ of affine $G$-subschemes of $V$ of type $\Gamma$, is
represented by a quasi-projective scheme. 
\end{theorem}

The invariant Hilbert scheme $\mathrm{Hilb}_\Gamma^G(V)$ is such a
quasi-projective scheme. 



\subsection{Closed points}

An algebraic variety endowed with an action of $G$ is said to be
\textit{spherical} 
if it is normal and it contains an open orbit of a Borel subgroup of $G$.

An affine $G$-variety $X$ is spherical if and only if its algebra of
regular functions $\mathbb C[X]$ is a \textit{multiplicity-free}
$G$-module, 
namely it is the direct sum of simple modules each of them occurring with
multiplicity one. 
In this case, let $\Gamma\subset\Lambda^+$ be the monoid such that
\begin{equation}
\mathbb C[X]\cong\bigoplus_{\lambda\in\Gamma}V(\lambda)^\ast
\end{equation}
as $G$-modules.
The normality of $X$ implies
$\mathbb Z \Gamma \cap \mathbb Q_{\leq 0} \Gamma = \Gamma$.

A spherical $G$-subvariety of $V$ with weight monoid equal to $\Gamma$
can be regarded as a closed point of $\mathrm{Hilb}_\Gamma^G(V)$.
It is said to be \textit{non-degenerate} if all its projections
into the isotypical components of $V$ are non-zero.

\begin{theorem}[\cite{AB} Corollary~1.17]
The non-degenerate spherical $G$-subvarieties of $V$ with weight monoid 
$\Gamma$ are parameterised by an open subscheme $\Hilb_\Gamma^G(V)_0$ of 
$\Hilb_\Gamma^G(V)$.
\end{theorem}

Let $X_0$ and $X_1$ be two quasi-affine spherical $G$-varieties,
$X_0$ is said to be a \textit{contraction}, or a 
\textit{degeneration} of $X_1$,
or vice versa that $X_1$ is a \textit{deformation} of $X_0$,
if $\C[X_0]$ is $G$-isomorphic to the graduate algebra associated to
a filtration of $\C[X_1]$.

Note that on a $G$-algebra $R\cong
\oplus_{\lambda\in\Gamma}V(\lambda)^\ast$ we have the maximal filtration 
$\mathfrak F_\mu=\oplus_{\lambda\leq\mu}V(\lambda)^\ast$, for all
$\mu\in\Gamma$, where $\lambda\leq\mu$ means that $\mu-\lambda$ 
is a sum of positive roots. 

\subsection{Toric action on the invariant Hilbert scheme}\label{Hilbert_action}

Let $Z(G)$ be the center of $G$. The invariant Hilbert scheme $\Hilb_\Gamma^G(V)$ is endowed with an action of the adjoint torus $T_\mathrm{ad}=T/Z(G)$ (\textit{see} \cite{AB}). The open subscheme $\Hilb_\Gamma^G(V)_0$ is stable for this action.
This action can be defined on multiplication laws of $R=\oplus_{\lambda\in\Gamma} V(\lambda)^\ast$ as follows. 

Let $X$ be a non-degenerate spherical $G$-subvariety of $V$ with weight monoid $\Gamma$. Its $G$-algebra of regular functions is isomorphic to $R$ defined above, as $G$-module. This gives a $G$-multiplication law in $R$, $m\colon R\otimes R\to R$. The map $m$ is direct sum of maps
$m_{\lambda,\mu}^\nu\colon V(\lambda)^\ast\otimes V(\mu)^\ast\to V(\nu)^\ast$
where $\lambda,\mu,\nu$ are in $\Gamma$, with $\nu\leq\lambda+\mu$.

For all $t\in T$, set
\begin{equation}\label{eq:action}
t.m=\bigoplus_{\lambda,\mu,\nu}t^{\lambda+\mu-\nu}m_{\lambda,\mu}^\nu.
\end{equation}

The limit of $t.m$ for $t^\alpha$ going to zero, for some $\alpha\in
S$, gives rise to a degeneration of $m$. If $t^\alpha$ goes to zero, for
all $\alpha\in S$, we obtain the so called \textit{Cartan multiplication},
that is $m_{\lambda,\mu}^\nu$ nonzero only if $\nu=\lambda+\mu$, 
corresponding to the maximal degeneration.

\subsection{Wonderful varieties}

An algebraic $G$-variety $X$ is said to be \textit{wonderful of rank $r$} if
\smallbreak
\noindent
{\rm(i)}\enspace
$X$ is smooth and complete
\smallbreak
\noindent
{\rm(ii)}\enspace
$G$ has an open $G$-orbit whose complement is the union of smooth prime $G$-divisors $D_i$ ($i=1,\ldots,r$) with normal crossings such that $\cap_1^r D_i\neq\emptyset$
\smallbreak
\noindent
{\rm(iii)}\enspace
if $x,x'$ are such that $\{i:x\in D_i\}=\{i:x'\in D_i\}$ then $G.x=G.x'$.
\smallbreak

A wonderful variety $X$ is projective and spherical (\emph{see}~\cite{L96}).

The \textit{spherical roots} of $X$ are the $T$-weights in the quotient
$T_z X/T_z Y$, 
where $Y$ is the unique closed $G$-orbit of $X$ and
$z\in Y$ is the unique point fixed by the opposite Borel subgroup of $B$.
The rank of $X$ is equal to the number of spherical roots of $X$.

Note that wonderful varieties of rank $0$ are homogeneous: they are the
flag varieties.

\section{Main results}

\subsection{The saturation assumption}

\begin{definition}[\cite{Pa}]
A submonoid of dominant weights, $\Gamma\subseteq\Lambda^+$, is said to be
\textit{saturated} if
\begin{equation}
\Z\Gamma\cap \Lambda^+=\Gamma.
\end{equation}
\end{definition}




The following characterisation of saturation for a free monoid is
straightforward. 

\begin{lemma}
\label{saturated}
Let $\lambda_1,\ldots,\lambda_s$ be linearly independent dominant weights.
The semigroup $\Gamma=\langle\lambda_1, \ldots, \lambda_s\rangle_\mathbb N$
is saturated if and only if
for every $i=1,\ldots,s$
there exists a  simple root $\alpha_{k_i}$ which is
not orthogonal only to the weight $\lambda_i$.
\end{lemma}


\subsection{The affine multicone over a flag variety}


Let $\lambda_1,\ldots,\lambda_s$ be dominant weights.
For all $i$, denote $P_{\lambda_i}$ the stabiliser
of the line $[v_{\lambda_i}]\in\P(V(\lambda_i))$
generated by a highest weight vector $v_{\lambda_i}\in V(\lambda_i)$.
Let $P=\cap_1^sP_{\lambda_i}$
that is, the stabiliser of $([v_{\lambda_1}],\ldots,[v_{\lambda_s}])$
in $\P(V(\lambda_1))\times\ldots\times\P(V(\lambda_s))$.
We have the following obvious inclusions
\begin{equation}
G/P\subset G/P_{\lambda_1}\times\ldots\times G/P_{\lambda_s}
\subset \P(V(\lambda_1))\times\ldots\times\P(V(\lambda_s)).
\end{equation}

Let $X_0$ denote the corresponding \textit{affine multicone} over the
partial flag variety $G/P$. It is the affine $G$-subvariety of
$V=\oplus_{i=1}^s V(\lambda_i)$ defined as 
\begin{equation}
X_0=\{v_1+\ldots+v_s\in V:
\exists\,([u_1],\ldots,[u_s])\!\in\! G/P \mbox{ with }v_i\!\in\![u_i]\,\forall i\}.
\end{equation}
It is a subvariety contained in the product $\mathcal
C_{\lambda_1}\times \ldots\times\mathcal C_{\lambda_s}$ of the affine cones
$\mathcal C_{\lambda_i}\subset V(\lambda_i)$ over $G/P_{\lambda_i}\subset
\P(V(\lambda_i))$. 

If the dominant weights $\lambda_1,\ldots,\lambda_s$ are supposed to be
linearly independent then 
the $G$-orbits of $X_0$ are parameterised by the subsets $I$ of
$\{1,\ldots,s\}$, 
\begin{equation}
G.v_{\underline\lambda}^I=\{v_1+\ldots+v_s\in X_0:v_i=0\,\forall i\in I\},
\end{equation}
where $v_{\underline\lambda}^I=\sum_{i\not\in I}v_{\lambda_i}$.
In particular, the $G$-orbit of
$v_{\underline\lambda}=(v_{\lambda_1},\ldots,v_{\lambda_s})$ is dense in
$X_0$. 

The $G$-variety $X_0$ is normal with an open $G$-orbit.
Further, its generic stabiliser is \textit{horospherical}, i.e.
it contains a maximal unipotent subgroup. 
In particular, the $G$-variety $X_0$ is spherical.

Denote $\Gamma$ the linear span of $\lambda_1, \ldots, \lambda_s$ over
$\mathbb N$ and, moreover, assume it to be saturated.
By Lemma\texttilde\ref{saturated} we get that, for any nonempty subset $I$,
the stabiliser $G_{v_{\underline\lambda}}$ is included in
$G_{v_{\underline\lambda}^I}$ 
with codimension at least two, therefore
\begin{equation}
\label{codimppty}
\mathrm{codim}_{X_0}(X_0\setminus G.v_{\underline\lambda})\geq2.
\end{equation}
The ring of regular functions $\C[X_0]$ is thus equal to
$\C[G/G_{v_{\underline\lambda}}]$, 
hence isomorphic to $\oplus_{\lambda\in\Gamma}V(\lambda)^\ast$ as a $G$-module.

The non-degenerate spherical $G$-subvariety $X_0$ of $V$
corresponds to the unique $T_\mathrm{ad}$-stable closed point in
$\mathrm{Hilb}_\Gamma^G(V)_0$
and any non-degenerate spherical $G$-subvariety of $V$ with weight monoid
$\Gamma$ is a deformation of $X_0$.


\subsection{Statements}

From now on, $\lambda_1, \ldots, \lambda_s$ are linearly independent
dominant weights generating a saturated monoid $\Gamma$ and $V$ is the
$G$-module given by $V(\lambda_1)\oplus\ldots\oplus V(\lambda_s)$. The
parabolic subgroup $P$ is the stabiliser of $\left([v_{\lambda_1}],\ldots,
  [v_{\lambda_s}]\right)$ in
$\P\left(V(\lambda_1)\right)\times\ldots\times\P\left(V(\lambda_1)\right)$, 
$X_0$ is the affine multicone over $G/P$ included in this multiprojective
space, as defined above.

\begin{theorem}
\label{sphericalweights}
As a $T_{\ad}$-module the tangent space of $\Hilb_\Gamma^G(V)$ at the point
$X_0$ is multiplicity free and  
its $T_{\ad}$-weights are spherical roots.
\end{theorem}

In the following we shall denote the above set of $T_{\ad}$-weights by
$\Sigma$. 

\begin{theorem}
\label{main}
\smallbreak
\noindent
{\rm (i)}\enspace
There exists a wonderful $G$-variety $X$ included in the multiprojective
space $\P\left( V(\lambda_1)\right) \times\ldots\times\P
\left(V(\lambda_s)\right)$ with set of spherical roots $\Sigma$
and with closed $G$-orbit $G/P$.
\smallbreak
\noindent
{\rm (ii)}\enspace
Let $\tilde{X}$ be the corresponding affine multicone over $X$.
Then there exists $v=(v_1,\ldots,v_s)\in \tilde{X}$, with
$([v_1],\ldots,[v_s])$ in the open $G$-orbit of $X$, such that
the subvariety $X_1=\overline{G.v}$ of $V$ is a deformation of $X_0$.
\end{theorem}

\begin{corollary}\label{affinespace}
The invariant Hilbert scheme $\Hilb_\Gamma^G(V)_0$
of non-degenerate spherical $G$-subvarieties of $V$
with weight monoid $\Gamma$
is an affine space.
Its dimension is equal to the rank of the wonderful variety $X$ of
Theorem~\ref{main}.
\end{corollary}

In particular, the dimension of the invariant Hilbert scheme is less or
equal than $s$, the rank of the free monoid $\Gamma$.
Moreover, the invariant Hilbert scheme is a point if and only if the
wonderful variety $X$ given above is of rank 0, namely $X=G/P$.

In Corollary~\ref{final} the cases where the invariant Hilbert scheme
$\Hilb_\Gamma^G(V)_0$ is not a point will be characterised by an explicit
condition on the monoid $\Gamma$. 

\begin{proof}
The spherical subvariety $X_1$ of $V$ given in Theorem~\ref{main} corresponds
to a closed point in $\Hilb_\Gamma^G(V)_0$. We shall prove that the
$T_{ad}$-orbit of $X_1$ in $\Hilb_\Gamma^G(V)_0$ is dense.

By Corollary~2.14 in \cite{AB} the normalisation of the closure of the
$T_\mathrm{ad}$-orbit of $X_1$ in $\Hilb^G_\Gamma(V)_0$ is isomorphic to a 
$T_\mathrm{ad}$-module whose $T_\mathrm{ad}$-weights are explicitly given. 
They are the linear independent generators of the free submonoid
\[
\mathsf M=\mathbb Q_{\geq0}\mathsf M'\cap \Z \mathsf M'\subseteq \mathbb N S,
\]
obtained from the submonoid
\[
\mathsf M'=\langle\lambda+\mu-\nu:\lambda,\mu,\nu 
\in\Gamma,\,m_{\lambda,\mu}^\nu\neq0\rangle_\mathbb N,
\]
where $m_{\lambda,\mu}^\nu\colon V(\lambda)^\ast\otimes V(\mu)^\ast\to 
V(\nu)^\ast$ and $m=\sum m_{\lambda,\mu}^\nu$ is the multiplication in 
$\C[X_1]$ (\emph{see} Paragraph~\ref{Hilbert_action}).

From Theorem~1.3 in \cite{K96} we know that such $T_{\ad}$-weights are 
exactly the spherical roots of the spherical varietx $X_1$, which are in
correspondece (but maybe not equal) to the spherical roots in $\Sigma$.
Therefore,  
by Theorem~\ref{sphericalweights}, the dimension of the orbit $T_{\ad}.X_1$
is equal to the dimension of the tangent space $T_{X_0}\Hilb_\Gamma^G(V)$. 
It follows that $\Hilb_\Gamma^G(V)_0$ is smooth,
hence equal to $\overline{T_{\ad}.X_1}$ and isomorphic to an affine space. 
\end{proof}

\begin{remark}
Let us keep the notation of the previous proof.
By Proposition~2.13 in \cite{AB} the monoid $\mathsf
M'$ is the weight monoid of 
$\overline{T_{\ad}.X_1}$ in $\Hilb_\Gamma^G(V)$, 
Therefore we have, a posteriori, that $\mathsf M'$ is
free.
\end{remark}

\begin{corollary}
\label{totalcone}
The union of all non-degenerate spherical $G$-subvarieties of $V$ with
weight monoid $\Gamma$ is equal to the affine multicone $\tilde X$
over the wonderful $G$-variety $X$ 
of Theorem~\ref{main}
included in 
$\P\left( V(\lambda_1)\right) \times\ldots\times\P \left( V(\lambda_s)\right)$.
\end{corollary}

The union of all non-degenerate spherical $G$-subvarieties of $V$ with
weight monoid $\Gamma$ is a closed subset of $V$.
By the same arguments as in the proof of Proposition~\ref{correspondence} 
one can get that for any point $([v'_1],\ldots,[v'_s])$ in $X$, 
$v'=(v'_1,\ldots,v'_s)$
is such that the subvariety $\overline{G.v'}$ of $V$ is a deformation of
$X_0$.
Then the above corollary will follow from Theorem~\ref{main}{\rm (i)} and
Corollary~\ref{affinespace}. 

\subsection{The case of a simple module}

As already mentioned in the introduction, S.~Jansou considers in \cite{J}
the case where the monoid $\Gamma$ is generated by a single dominant weight
$\lambda$ and proves that the corresponding Hilbert scheme is a point or an
affine line. Moreover, he discovers a correspondence with the list of one
dimensional Hilbert schemes and the list of two-orbit varieties whose
closed orbit is an ample divisor.


The two-orbit varieties $X'$ which appear in this setting are (by
definition) wonderful of rank one.

The union of the equivariant
deformations of $X_0$ (here $X_0$ is the cone of the primitive vectors in
$V(\lambda)$) is obtained by taking the affine cone over the wonderful
variety $X$ (still of rank one) whose generic stabiliser is the normaliser
of the generic stabiliser of the former rank one wonderful variety $X'$.

\subsection{The model case}

On the other side D.~Luna in~\cite{L05} considers the case where the monoid
$\Gamma$ is equal to $\Lambda^+$. 
More precisely he is concerned with the so
called \textit{model homogenous spaces}, that are the quasi-affine
(spherical) $G$-homogeneous spaces whose ring of regular functions
is isomorphic to $\oplus_{\lambda\in\Lambda^+}V(\lambda)^\ast$ as a $G$-module. 
He obtains that their isomorphism classes are
parameterised by the $G$-orbits of a certain wonderful variety for a given connected reductive group $G$. 
Further, this wonderful variety can be constructed explicitly case by case.

\subsection{On the saturation assumption}\label{par:onsat}

Let $H$ be a spherical subgroup of $G$, 
that is the homogeneous space $G/H$ has an open $B$-orbit.
The \textit{horospherical contraction} $\widehat H$ of $H$ (\cite{Pa},
\cite{K90}), in the case of a quasi-affine $G/H$, is defined as follows. 

Let $\Gamma(G/H)$ be the submonoid of highest weights in $\C[G/H]$,
then $\widehat H$ is the intersection of the stabilisers $G_{v_\lambda}$ 
for all weights $\lambda$ in $\Gamma(G/H)$.

The subgroup $\widehat H$ is horospherical. 
Moreover, $\dim G/\widehat H=\dim G/H$ and (\cite{Pa}
Proposition\texttilde1.5) 
\begin{equation}
\Gamma(G/\widehat H)=\Z\Gamma(G/H)\cap\Lambda^+.
\end{equation}
Therefore, if $\Gamma$ is saturated, 
the $G$-algebra $\C[G/\widehat H]$ is isomorphic to the graduate algebra 
associated to the maximal filtration $\mathfrak F_\mu$.
Hence $\mathrm{Spec}\,\C[G/H]$ is a deformation of 
$\mathrm{Spec}\,\C[G/\widehat H]$.

This can be applied in our setting. Indeed, if $\Gamma$ is free and saturated and
$X_0$ is the affine multicone defined above, then
a spherical $G$-subvariety $X_1$ of $V$ is a deformation of
$X_0$ if and only if its open $G$-orbit is a deformation of
$G.v_{\underline\lambda}$ (\emph{see} \cite{BK}).

\section{The tangent space of the invariant Hilbert
  scheme}\label{sct:tangentspace}

The purpose of this section is the proof of Theorem~\ref{sphericalweights}.
We first compute the $T_{\ad}$-weights of the tangent space of
$\Hilb_\Gamma^G(V)$ at the multicone $X_0$
$\mathrm T_{X_0}\Hilb_\Gamma^G(V)$ and then prove its multiplicity freeness.

Recall that $\lambda_1,\ldots,\lambda_s$ are linearly independent dominant
weights which generate a saturated submonoid $\Gamma$ and that
$V$ is the $G$-module $V(\lambda_1)\oplus\ldots\oplus V(\lambda_s)$.

Since $\mathrm{codim}_{X_0}(X_0\setminus G.v_{\underline\lambda})\geq2$
(\emph{see}~(\ref{codimppty})) we can use the following characterisation 
of the tangent space of the invariant Hilbert scheme at the multicone $X_0$
given by Proposition 1.15 of~\cite{AB} 
\begin{equation}
\mathrm T_{X_0}\Hilb_\Gamma^G(V)\cong
(V/\mathfrak{g}. v_{\underline{\lambda}})^{G_{v_{\underline\lambda}}}
\end{equation}
where $\mathfrak g$ denotes the Lie algebra of $G$ and
$v_{\underline{\lambda}}=v_{\lambda_1}+\ldots+v_{\lambda_s}$,
with $v_{\lambda_i}$ a highest weight vector in $V(\lambda_i)$.

\subsection{Action of the adjoint torus - its weights}


Consider the following action of the adjoint torus $T_{\ad}=T/Z(G)$ on $V$
defined on every $T$-weight space $V(\lambda_i)_\mu$ by
\begin{equation}
t.v=\lambda_i(t)\mu(t)^{-1}v
\end{equation}
(\emph{see} Section~2.1 in \cite{AB}). It yields a $T_{\ad}$-action
on the vector space
$(V/\mathfrak g.v_{\underline\lambda})^{G_{v_{\underline\lambda}}}$.
This action is often called the \emph{normalised action} of $T_{\ad}$ on
$(V/\mathfrak g.v_{\underline\lambda})^{G_{v_{\underline\lambda}}}$.
It corresponds to the differential of the $T_\mathrm{ad}$-action
on the invariant Hilbert scheme $\Hilb_\Gamma^G(V)$
defined in (\ref{eq:action}) (\textit{see loc.\ cit.}).

The stabiliser $G_{v_{\underline\lambda}}$ is generated by its maximal torus,
intersection of the kernels of the characters $\lambda_i$,
and its unipotent part.
The $T_{\ad}$-module
$(V/\mathfrak g.v_{\underline\lambda})^{G_{v_{\underline\lambda}}}$
is thus equal to
\begin{equation}
(V/\mathfrak g.v_{\underline\lambda})^{\mathfrak g_{v_{\underline\lambda}}}
_{\langle\lambda_1,\ldots,\lambda_s\rangle},
\end{equation}
that is, the subspace spanned by the $T_{\ad}$-semi-invariant classes in
$(V/\mathfrak g.v_{\underline\lambda})^{\mathfrak g_{v_{\underline\lambda}}}$
whose $T_\mathrm{ad}$-weights are integral linear combinations of the weights
$\lambda_i$.

Note that representatives of $T_{\ad}$-weight vectors of
$(V/\mathfrak g.v_{\underline\lambda})^{\mathfrak g_{v_{\underline\lambda}}}$
with weight $\gamma$ can be chosen in
$\oplus_i V(\lambda_i)_{\lambda_i-\gamma}$.

\begin{theorem}
\label{sphericalweightslist}
Let $\lambda_1,\ldots,\lambda_s$ be linearly independent dominant weights
generating a saturated cone $\Gamma$.
Consider the $G$-module $V=\oplus_1^sV(\lambda_i)$.
The $T_\mathrm{ad}$-weights of $(V/\mathfrak{g}. v_{\underline\lambda})^{G_{v_{\underline\lambda}}}$ for the normalised action belong to Table~\ref{tab:weights}.
\end{theorem}

\begin{howtoreadtableone}
A $T_\mathrm{ad}$-weight $\gamma$ is an element of the root lattice.
In the table it is written as linear combination
of simple roots.
The simple roots are labelled according to the type of the support
of the $T_\mathrm{ad}$-weight $\gamma$.
\end{howtoreadtableone}

\begin{table}[tbh]
\caption{$T_\mathrm{ad}$-weights in $(V/\mathfrak g.v_{\underline\lambda})^{G_{v_{\underline\lambda}}}$.}
\label{tab:weights}
\begin{center}
\begin{tabular}{|l|l|}
\hline
Type of support & \multicolumn{1}{|c|}{Weight} \\
\hline
\hline
$\mathsf A_1\times\mathsf A_1$ & $\alpha+\alpha'$ \\
\hline
$\mathsf A_n$ & $\alpha_1+\ldots+\alpha_n$, $n\geq2$ \\
& $2\alpha$, $n=1$ \\
& $\alpha_1+2\alpha_2+\alpha_3$, $n=3$ \\
\hline
$\mathsf B_n$, $n\geq2$ & $\alpha_1+\ldots+\alpha_n$ \\
& $2\alpha_1+\ldots+2\alpha_n$ \\
& $\alpha_1+2\alpha_2+3\alpha_3$, $n=3$ \\
\hline
$\mathsf C_n$, $n\geq3$ & $\alpha_1+2\alpha_2+\ldots+2\alpha_{n-1}+\alpha_n$ \\
\hline
$\mathsf D_n$, $n\geq4$ & $2\alpha_1+\ldots+2\alpha_{n-2}+\alpha_{n-1}+\alpha_n$ \\
& $\alpha_1+2\alpha_2+2\alpha_3+\alpha_4$, $n=4$ \\
& $\alpha_1+2\alpha_2+\alpha_3+2\alpha_4$, $n=4$ \\
\hline
$\mathsf F_4$ & $\alpha_1+2\alpha_2+3\alpha_3+2\alpha_4$ \\
\hline
$\mathsf G_2$ & $4\alpha_1+2\alpha_2$ \\
& $\alpha_1+\alpha_2$ \\
\hline
\end{tabular}
\end{center}
\end{table}

\begin{corollary}
The $T_\mathrm{ad}$-weights of $(V/\mathfrak{g}. v_{\underline\lambda})^{G_{v_{\underline\lambda}}}$ are spherical roots.
\end{corollary}

\begin{proof}
From the list stated in Table~\ref{tab:weights},
one checks easily that the elements of $\Sigma$ are indeed spherical roots
(compare with Table~1 in~\cite{W} which lists all the spherical roots).
\end{proof}

The following is devoted to the proof of Theorem~\ref{sphericalweightslist}.


\begin{lemma}\label{pairlemma}
Let $[v]\neq0$ be a $T_{\ad}$-weight vector
in $(V/\mathfrak{g}. v_{\underline\lambda})^{G_{v_{\underline\lambda}}}$.
Then there exists a pair $(\alpha,\beta)$
of positive roots with $\alpha$ simple such that (up to a rescaling of $v$)
\begin{equation}
\label{pair}
X_\alpha. v\neq 0\quad \mbox{and}\quad X_\alpha. v=X_{-\beta}. v_{\underline\lambda}.
\end{equation}

\end{lemma}

\begin{proof}
The class $[v]$ being different from $0$,
the vector $v$ is not a linear combination of highest weight vectors
$v_{\lambda_i}$. 
Hence, there exists $\alpha\in S$ such that $X_\alpha.v\neq 0$.
Since $\mathfrak g_{v_{\underline\lambda}}.v\subset\mathfrak g.v_{\underline\lambda}$,
we have $X_\alpha.v\in\mathfrak g.v_{\underline\lambda}$.

Assume $X_\alpha.v$ is a linear combination of highest weight vectors.
The simple root $\alpha$ is thus the $T_{\ad}$-weight of $[v]$
and $v$ can be chosen in $\oplus_iV(\lambda_i)_{\lambda_i-\alpha}$.
Hence, there exists $j$ such that $(\lambda_j,\alpha)\neq0$.
Because of saturation, by Lemma\texttilde\ref{saturated},
such a weight $\lambda_j$ is unique,
since $\alpha$ is an integral combination of the weights $\lambda_i$.
This means that the $i$-th component of $v$ in $V(\lambda_i)$ is zero,
for all $i\neq j$.
Hence, $v\in\C X_{-\alpha}.v_{\underline\lambda}$
which is not possible because $[v]\neq 0$.
Therefore, there exists $\beta\in\Phi^+$ such that
$X_\alpha.v\in\C X_{-\beta}.v_{\underline\lambda}$.
\end{proof}

\begin{remark}
\label{weight}
The sum $\alpha+\beta$ of any pair $(\alpha,\beta)$ of positive roots
(with $\alpha$ not necessarily simple)
which satisfies condition\texttilde(\ref{pair}) is equal to the weight of
$[v]$. 
\end{remark}

\begin{proposition}
\label{insidesupp}
Let $\gamma$ be a $T_{\ad}$-weight
of some $(V/\mathfrak g .v_{\underline\lambda})^{G_{v_{\underline\lambda}}}$.
Suppose there exists a simple root $\delta$ in the support of $\gamma$
such that $\gamma-\delta\not\in\Phi$.
Then $(\gamma,\delta)\geq 0$.
Moreover if $(\gamma,\delta)=0$ then $(\lambda_i,\delta)=0$
for all $i$.
\end{proposition}

\begin{proof}
Let us prove the two assertions together.

We shall proceed by contradiction.
Suppose the root $\delta$ is such that $(\gamma,\delta)\leq0$
and there exists $j$ such that $(\lambda_j,\delta)\neq 0$
(that is obviously true if $(\gamma,\delta)\neq 0$).

By Lemma\texttilde\ref{saturated} there exist a simple root $\alpha$ and
a weight $\lambda_k$ with $(\gamma,\alpha)>0$,
$(\lambda_k,\alpha)>0$ and $(\lambda_i,\alpha)=0$ for all $i\neq k$.
Since $\lambda_k$ occurs with a positive coefficient in the writing of $\gamma$
as an integral linear combination of the weights $\lambda_i$,
if $(\lambda_k,\delta)$ is nonzero
$(\lambda_i,\delta)$ can not be zero for all $i\neq k$.
Hence, we can suppose $(\lambda_j,\alpha)=0$.

Let $[v]\in(V/\mathfrak g .v_{\underline\lambda})^{G_{v_{\underline\lambda}}}$ be nonzero
of $T_{\ad}$-weight $\gamma$.
Consider $v_j\in V(\lambda_j)$ the $j$-th component
of some representative $v$ of $[v]$.
Note that we have $(\lambda_j-\gamma,\alpha)<0$.
Hence, if $v_j\neq 0$ then $X_{\alpha}. v_j\neq 0$,
identity~(\ref{pair}) holds with the pair $(\alpha,\gamma-\alpha)$
and in particular $X_\alpha.v_j=X_{-\gamma+\alpha}.v_{\lambda_j}$.
This implies that $v_j$ is proportional to
$X_{-\alpha}X_{-\gamma+\alpha}.v_{\lambda_j}$
which equals $X_{-\gamma}. v_{\lambda_j}$ in case
$\gamma\in\Phi$ and $0$ otherwise.
In any case, we thus have: $v_j\in\mathfrak g.v_{\lambda_j}$
hence we can choose $v'\in [v]$ such that $v'_j=0$.

Take $\alpha'$ and $\beta'$ that fulfill (\ref{pair}) for the vector $v'$,
we have: $(\lambda_j,\gamma-\alpha')=0$.
Since $\delta\in\supp\gamma$ and $(\lambda_j,\delta)\neq 0$,
the only possibility is that $\delta=\alpha'$.
But $\gamma-\delta$ is not a root,
hence (\ref{pair}) can not hold with $\delta$.
This yields the contradiction.
\end{proof}

\begin{lemma}\label{restriction}
Let $\gamma$ be a $T_{\ad}$-weight of 
$(V/\mathfrak g.v_{\underline\lambda})^{G_{v_{\underline\lambda}}}$
and $L$ be the standard Levi subgroup $L$ 
associated to the set of simple roots $\supp\gamma$. 
Denote $W$ the $L$-submodule of $V$ generated by $v_{\pi{\underline\lambda}}$.
Then 
\begin{equation}
(V/\mathfrak g.v_{\underline\lambda})^{G_{v_{\underline\lambda}}}_\gamma
\cong(W/\mathfrak l.v_{\pi(\underline\lambda)})
^{L_{v_{\pi(\underline\lambda)}}}_\gamma .
\end{equation}
as $T_{\ad}$-modules, where $\pi$ is the projection of $\Lambda$ onto the
weight lattice of $L$.
\end{lemma}

\begin{proof}
The isomorphism of the lemma is given by sending 
$[v]\in V/\mathfrak g.v_{\underline\lambda}$ to the class of $v$ in  
$W/\mathfrak l.v_{\underline\lambda}$
where $v\in \oplus_i V(\lambda_i)_{\lambda_i -\gamma}$. 
\end{proof}


\begin{lemma}\label{ruleout}
Let $\gamma$ be a $T_{\ad}$-weight of some $(V/\mathfrak g .v_{\underline\lambda})^{G_{v_{\underline\lambda}}}$.
Suppose $\gamma$ is not a root.
Consider $\alpha$ and $\beta$ two positive roots which fulfill condition~(\ref{pair}).
If $\beta=\beta_1+\beta_2$ with $\beta_1$ and $\beta_2$ positive roots
then $\alpha+\beta_1$ or $\alpha+\beta_2$ is a root.
\end{lemma}

\begin{proof}
Let $v\in V$ satisfy (\ref{pair}) with the given roots $\alpha$ and $\beta$,
in particular $\gamma=\alpha+\beta$.
Assume $\alpha+\beta-\beta_i$ is not a root for $i=1,2$.
Since $X_{\beta_i}.v\in\mathfrak g.v_{\underline\lambda}$ we have: $X_{\beta_i}.v=0$ for $i=1,2$ hence $X_\beta.v=0$.
It follows along with $\gamma\not\in\Phi$ and the property~(\ref{pair}):
$0=X_\alpha X_\beta. v=X_\beta X_\alpha. v=X_\beta X_{-\beta}. v_{\underline\lambda}$.
In particular, we have: $X_{-\beta}.v_{\underline\lambda}=0$ which contradicts $0\neq X_\alpha.v=X_{-\beta}.v_{\underline\lambda}$
given by~(\ref{pair}).
\end{proof}


\begin{lemma}\label{simple1}
Let $[v]$ be a nonzero semi-invariant class in $(V/\mathfrak
g.v_{\underline\lambda})^G_{v_{\underline\lambda}}$ 
and $(\alpha,\beta)$ pair of positive roots, with $\alpha$ simple,
fulfilling condition\texttilde\eqref{pair} for the vector $v$.
If $\alpha$ and $\beta$ are orthogonal and $\alpha+\beta$ is not a root,
namely the pair $(\alpha,\beta)$ generates a root subsystem of type
$\mathsf A_1\times\mathsf A_1$, 
then there exists a unique $j$ such that $\lambda_j$ is not orthogonal to
$\alpha+\beta$. 
\end{lemma}

\begin{proof}
Consider the root subsystem $\Phi'$ generated by $\alpha$ and $\beta$
and apply Lemma\texttilde\ref{restriction}.
\end{proof}

\begin{lemma}\label{simple2}
Let $\gamma$ be the weight of a nonzero class $[v]\in(V/\mathfrak g.v_{\underline\lambda})^G_{v_{\underline\lambda}}$.
If one of the following conditions is fulfilled, there exists
a unique $j$ such that $\lambda_j$ is not orthogonal to $\gamma$.
\smallbreak
\noindent
{\rm (i)}\enspace
There exists a unique root $\delta_0\in\supp\gamma$ such that $(\gamma,\delta_0)\neq0$
and there exist no root $\delta\in\supp\gamma$ such that $(\gamma,\delta)=0$ and $\gamma-\delta\in\Phi$.
\smallbreak
\noindent
{\rm(ii)}\enspace
There exists a unique simple root $\alpha$ such that
$(\gamma,\alpha)>0$ and $\gamma-\alpha\in\Phi$.
Moreover, for all $\delta\in\supp\gamma\setminus\{\alpha\}$
we have:
$(\gamma,\delta)\leq0$ and $\gamma-\delta\not\in\Phi$.
\end{lemma}

\begin{proof}
By saturation and Proposition~\ref{insidesupp}.
\end{proof}

\paragraph{Proof of Theorem~\ref{sphericalweightslist}.}
%
Pick a $T$-weight $\gamma$ of $(V/\mathfrak{g}. v_{\underline\lambda})^{G_{v_{\underline\lambda}}}$.
We know by Lemma~\ref{pairlemma} that the weight $\gamma$
is a sum of two positive roots,
say $\alpha$ and $\beta$ with $\alpha$ simple.
If $\alpha$ is orthogonal to the support of $\beta$,
then $\beta$ is also simple by Lemma~\ref{ruleout}.

This allows to reduce to the case of a simple group $G$.

Finally most of the choices for $\gamma=\alpha+\beta$ can be ruled out:
in some cases we can use Lemma~\ref{ruleout}, in other cases it is possible
to prove that there is only one dominant weight $\lambda_i$ not orthogonal
to $\gamma$ (applying Lemma~\ref{simple1} and Lemma~\ref{simple2}) hence we
can use Lemma~\ref{restriction} and Theorem~1.1 in \cite{J}. 

We end up with the list of Table~\ref{tab:weights}. 

\subsection{Multiplicity freeness}

\begin{theorem}
Each $T_{\ad}$-weight space of 
$(V/\mathfrak g.v_{\underline{\lambda}})^{G_{v_{\underline{\lambda}}}}$
is one dimensional.
\end{theorem}

We first work out the case of irreducible modules $V$ then in combination with
Lemma~\ref{restriction} we get to the general case.

\begin{proposition}
\label{simplecase}
Suppose $V$ is an irreducible $G$-module, 
namely $V=V(\lambda)$ for a dominant weight $\lambda$.
Then the $T_{\ad}$-weight spaces of  
$(V/\mathfrak g.v_{\underline\lambda})^{G_{v_{\underline\lambda}}}$
are one dimensional.
\end{proposition}
 
The previous statement was already proved in \cite{J}. 
We present here a more concise proof.

\begin{proof} 
From Table~\ref{tab:weights}, we know that 
if the weight $\gamma$ is not a root 
then there exists exactly one simple root $\delta$ 
such that $(\gamma,\delta)\neq 0$. 
Hence, by Proposition~\ref{insidesupp}, 
the dominant weight $\lambda$ is (up to a scalar) 
the fundamental weight associated to $\delta$.
Moreover, for such a root $\delta$, we can check that 
the coefficient corresponding to $\delta$ 
in the writing of $\gamma$ is $2$ or $3$. 
Therefore, if this coefficient is $2$, 
the following vectors are clearly basis vectors 
of the weight space $V_{\gamma-\lambda}$:
$$
X_{-\gamma+\nu}X_{-\nu}.v_\lambda,\enspace 
$$ 
where $\{\nu,\gamma-\nu\}$ 
runs over the $2$-set 
of positive roots such that $\nu$ and $\gamma-\nu$ 
contain $\delta$ in their support.
This holds also in case the coefficient of $\delta$ is $3$.
Indeed, we have then $\gamma=\alpha_1+2\alpha_2+3\alpha_3$ 
of type $\mathsf B_3$.
Consider the vectors
$X_{-\gamma+\nu_1+\nu_2}X_{-\nu_1}X_{-\nu_2}.v_\lambda$ 
with $\{\nu_1,\nu_2,\gamma-(\nu_1+\nu_2)\}$ 
running over $3$-sets of positive roots 
$\nu_1$, $\nu_2$ and $\gamma-(\nu_1+\nu_2)$ 
which all contain $\delta$ in their support.
Then it is easy to see that these vectors 
are linear combinations of the above vectors 
$X_{-\gamma+\nu}X_{-\nu}.v_\lambda$.

Applying now the equalities $X_\nu.v=0$ for all 
$\nu\in\supp\gamma\setminus\{\delta\}$, we obtain the proposition.

If $\gamma$ is a root then we are left with $\gamma=\alpha_1+\ldots+\alpha_n$, 
$\gamma=\alpha_1+2(\alpha_2\ldots)+\alpha_n$ or 
$\gamma=\alpha_1+2\alpha_2+3\alpha_3+2\alpha_4$
respectively of type $\mathsf A_n$, $\mathsf C_n$ and $\mathsf F_4$.
Indeed, type $\mathsf B_n$ (resp. $\mathsf G_2$) is easily ruled out 
since the weight space $V(\lambda)_{\lambda-\gamma}$ is spanned by 
$X_{-\gamma}.v_{\lambda}$
(resp. $\gamma=\omega_2-\omega_1$ is not up to a scalar a dominant weight).
For the possible remaining weights $\gamma$, 
we take as  basis vectors of $V(\lambda)_{\lambda-\gamma}$
the same vectors $X_{-\gamma+\nu}X_{-\nu}.v_\lambda$ as above 
($\delta=\alpha_1,\alpha_2$ or $\alpha_4$ respectively) 
along with $X_{-\gamma}.v_\lambda$
and conclude like before.
\end{proof}

\begin{proof}
Let $\gamma$ be a $T_{\ad}$-weight of 
$(V/\mathfrak g.v_{\underline{\lambda}})^{G_{v_{\underline{\lambda}}}}$.

We first consider the case where $G$ is of type $\mathsf G_2$ 
and $\gamma=\alpha_1+\alpha_2$, 
we find as single possibility for $\underline\lambda$:
$\lambda_1=\omega_1$ and $\lambda_2=\omega_2$ 
and the tangent space is one dimensional.

Suppose now that if $\gamma= \alpha_1+\alpha_2$ 
then its support is not of type $\mathsf G_2$.
From Table~\ref{tab:weights}, we get:
$$
(\gamma,\alpha)\geq 0 \mbox{ for every }\alpha\in \supp\gamma .
$$
Note that this can be proved also by means of Proposition~\ref{insidesupp}.

One can check in Table~\ref{tab:weights} 
that there are either one or two simple roots $\delta$
such that $(\gamma,\delta)>0$. Such roots are necessarily in $\supp\gamma$.

Consider first the case where there exists a unique simple root, say
$\delta$, such that $(\gamma,\delta)>0$.
Hence, because of saturation, there exists a unique dominant weight 
$\lambda_i$ such that $(\lambda_i,\delta)\neq 0$.
By assumption we have: $(\gamma,\alpha)=0$ for every 
$\alpha\in\supp\gamma\setminus\{\delta\}$.

If $(\underline{\lambda},\alpha)=0$ 
for every simple root $\alpha\in\supp\gamma\setminus\{\delta\}$,
we have: 
$$
\oplus_j V(\lambda_j)_{\lambda_j-\gamma}=V(\lambda_i)_{\lambda_i-\gamma}.
$$
and the theorem follows from the reduction to the Levi subgroup 
associated to the support of $\gamma$ and Proposition~\ref{simplecase}.

Suppose that $(\underline{\lambda},\alpha)\neq 0$ 
for a simple root $\alpha\in\supp\gamma\setminus\{\delta\}$.
By Proposition~\ref{insidesupp}, 
the root $\alpha$ is such that $\gamma-\alpha\in\Phi$.
Note that the weight $\gamma$ is obviously a root 
and the root $\alpha$ is unique.
By Table~\ref{tab:weights}, we have even more precisely 
that $\gamma=\alpha_1+\ldots+\alpha_n$, 
$\gamma=\alpha_1+2(\alpha_2+\ldots)+\alpha_n$ or
$\gamma=\alpha_1+2\alpha_2+3\alpha_3+2\alpha_4$ 
respectively of type $\mathsf B_n$, $\mathsf C_n$ and $\mathsf F_4$.
The root $\alpha$ is respectively $\alpha_n$, $\alpha_1$ and $\alpha_3$.
Then  
$\underline{\lambda}$ have one of the following shapes 
(after restriction to the Levi subgroup associated to the support of $\gamma$):
$$
(\omega_\delta,a_1\omega_\alpha,\ldots,a_r\omega_\alpha) 
\enspace\mbox { or }\enspace 
(\omega_\delta+b\omega_\alpha,b_1\omega_\alpha,\ldots,b_r\omega_\alpha) 
$$
where the $a_i$'s, $b$ and the $b_j$'s are positive integers and $r\geq 1$.

Note that the case of type $\mathsf F_4$ is special. 
Only the first shape of $\underline{\lambda}$ is possible a priori, 
with $r=1$, but in this case also the tangent space is $0$-dimensional.
We postpone it in a separated lemma below.

In type $\mathsf B_n$ and $\mathsf C_n$ the weight spaces 
$V(\lambda_i)_{\lambda_i-\gamma}$ are clearly spanned by 
$X_{-\gamma}.v_{\lambda_i}$ for any positive integer 
$a_i$, $b_i$ ($i=1,\ldots,r$).
Let $[v]$ be a $T_{\ad}$-vector of weight $\gamma$ with 
$v\in \oplus_i V(\lambda_i)_{\lambda_i -\gamma}$.
We can thus suppose all its components  but the first one to be $0$. 
Indeed, one can take as new representative for $[v]$, 
the vector $v-kX_{-\gamma+\delta}.v_{\underline\lambda}$ 
if one of the $v_j$ ($j>1$) is such that
$0\neq v_j=kX_{-\gamma+\delta}.v_{\lambda_j}$ ($k$ being a non-zero scalar).

If $\lambda_1=\omega_\delta$ with $\gamma=\alpha_1+\ldots+\alpha_n$ 
(here $\delta=\alpha_1$)
then the tangent space is clearly one dimensional 
and we can take as $T_{\ad}$-weight vector
the class of $(0,\ldots,0,X_{-\gamma}.v_{\underline\lambda})$.
Remark that only in this case the weight space 
$V(\lambda_1)_{\lambda_1-\gamma}$ is one dimensional. 
For the other cases under consideration, 
one basis of $V(\lambda_1)_{\lambda_1-\gamma}$ 
is given by the following vectors:
$$
X_{-\gamma}. v_{\lambda_1}, X_{-\gamma+\nu} X_{-\nu}.v_{\lambda_1}
$$
with $\{\delta\}\subset\supp\nu\subset\supp\gamma\setminus\{\alpha\}$.

Recall that $X_\nu.v=0$ for all $\nu\in\supp\gamma\setminus\{\alpha,\delta\}$ 
since $\gamma-\nu\not\in\Phi$.
Writing $v$ as a linear combination of the above basis vectors 
and using these equalities, 
we obtain that the tangent space is multiplicity free.

Consider now the case where there exist two distinct simple roots, 
say $\delta_1$ and $\delta_2$,
such that $(\gamma,\delta_1)>0$ and $(\gamma,\delta_2)>0$.
According to Table~\ref{tab:weights}, this assumption is only valid 
for $\gamma=\alpha_1+\ldots+\alpha_n$ of type $A_n$.
After reduction to the Levi subgroup associated to $\supp \gamma$, 
we get as possibilities for $\underline\lambda$:
$$
\omega_1+\omega_n,\enspace (\omega_1,\omega_n),\enspace 
(a\omega_1+\omega_n,b\omega_1),
\enspace (a\omega_1+\omega_n,b\omega_1,c\omega_1)
$$
where $a,b,c$ are positive integers. 
Note that the last case can only occur in type $\mathsf D$ and $\mathsf E$.
In all these cases, one shows similarly as before that the components 
in $V(\lambda_i)$, $i=2,3$, of a $T_{\ad}$-weight vector can be taken 
to be equal to $0$.
Then straightforward computations lead to the proposition.
\end{proof}

\begin{lemma}\label{f4}
If $G$ is of type $\mathsf F_4$, 
$\lambda_1=\omega_4$ and $\lambda_2=a\omega_3$ with $a>0$, then the space 
$(V/\mathfrak g.v_{\underline{\lambda}})^{G_{v_{\underline{\lambda}}}}$
is trivial.
\end{lemma}

\begin{proof}
We can take $v_2=0$ since it should be 
in $(V(\lambda_2)/\mathfrak g.v_{\lambda_2})^{G_{v_{{\underline\lambda}}}}$ 
and in fact in $(V(\lambda_2)/\mathfrak g.v_{\lambda_2})^{G_{v_{\lambda_2}}}$ 
which is trivial as proved previously. 
We have also 
$v_1\in (V(\lambda_1)/\mathfrak g.v_{\lambda_1})^{G_{v_{\lambda_1}}}$ 
and recall that in this case the space is one dimensional.
Suppose $v_1\neq 0$ then since $X_{\alpha_3}.v_1\neq 0$, it should be 
(up to a scalar) equal to $X_{-\gamma+\alpha_3}.v_{\lambda_1}$.
Hence $X_{\alpha_3}.v_2$ (which is $0$) 
is not equal to $X_{-\gamma+\alpha_3}.v_{\lambda_2}$. 
\end{proof}

\section{Spherical systems}

After having recalled the axiomatic definition of a spherical system given in
~\cite{L01}, 
we show how we can associate naturally to the given monoid $\Gamma$ a
spherical system. 
Further, we prove in the last paragraph that $\Gamma$ fulfills a certain
regularity property.

\subsection{Definition of spherical systems}
\label{sphsys}

\begin{definition}
A \textit{spherical system} for $G$
is a triple consisting of a subset $S^p$ of simple roots,
a set $\Sigma$ of spherical roots for $G$
(namely $T$-weights that are the spherical root of a rank one wonderful $G$-variety) and
a set $\mathbf A$ endowed with a map $\rho\colon\mathbf A\to\Xi^\ast$, where $\Xi=\Z\Sigma$,
which satisfies the following properties.
\smallbreak\noindent(A1)\enspace For every $D \in \mathbf A$ and $\gamma \in \Sigma$
we have $\langle \rho (D),\gamma \rangle \leq 1$, and if $\langle\rho(D),
\gamma\rangle=1$ then $\gamma \in S\cap\Sigma$.
\smallbreak\noindent(A2)\enspace For every $\alpha \in \Sigma \cap S$,
$\mathbf A$ contains two elements, $D_\alpha^+$ and $D_\alpha^-$,
such that $\langle\rho(D_\alpha^\pm),\alpha\rangle=1$.
Moreover
$\langle \rho (D_\alpha^+),\gamma \rangle + \langle \rho (D_\alpha^-),\gamma \rangle =
(\gamma,\alpha^\vee)$, for every $\gamma \in \Sigma$.
\smallbreak\noindent(A3)\enspace The set $\mathbf A$ is the union of $\{D_\alpha^+,D_\alpha^-\}$
for all $\alpha\in\Sigma \cap S$.
\smallbreak\noindent($\Sigma 1$)\enspace If $2\alpha \in \Sigma \cap 2S$ then $\frac{1}{2}\langle\alpha^\vee, \gamma \rangle$
is a nonpositive integer for every $\gamma \in \Sigma \setminus \{ 2\alpha \}$.
\smallbreak\noindent($\Sigma 2$)\enspace If $\alpha, \beta \in S$ are orthogonal and $\alpha + \beta \in \Sigma$
or $\frac{1}{2}(\alpha + \beta) \in \Sigma$
then $(\gamma,\alpha ^\vee) = (\gamma,\beta ^\vee)$
for every $\gamma \in \Sigma$.
\smallbreak\noindent(S)\enspace For every $\gamma \in \Sigma$, there exists a rank one wonderful $G$-variety $X$
with $\gamma$ as spherical root and $S^p$ equal to the set of simple roots associated to $P_X$.
\end{definition}

\begin{remark}\label{combinatorialaxiom}
The last axiom (S) of spherical system is equivalent, 
in type different from $\mathsf F$, to the following.
For all $\gamma\in\Sigma$,
\begin{equation}
\{\alpha\in S:(\gamma,\alpha)=0\ \mbox{and}\ \gamma-\alpha\not\in\Phi\}\cap
\supp\gamma\subseteq S^p\subseteq\{\alpha\in S:(\gamma,\alpha)=0\}.
\end{equation}
\end{remark}

\subsection{The spherical system attached to $\Gamma$}
\label{par:sphericalsystem}

Let $\Gamma$ be a saturated cone generated by linearly independent
dominant weights $\lambda_1,\ldots,\lambda_s$.
Let $V$ be the finite dimensional $G$-module $\oplus_iV(\lambda_i)$.
We denote by $S^p$ be the subset of simple roots orthogonal to $\lambda_i$ for all $i$
and by  $\Sigma$ be the set of $T_\mathrm{ad}$-weights
in $(V/\mathfrak g._{v_{\underline\lambda}})^{G_{v_{\underline\lambda}}}$.

\begin{theorem}\label{sphericalsystem}
The triple ($S^p$, $\Sigma$, $\emptyset$) is a spherical system for $G$.
\end{theorem}

The two following results are directly involved in the proof of
Theorem\texttilde\ref{sphericalsystem}
and Proposition\texttilde\ref{char}.
The second one is a consequence of the first one
and of Theorem\texttilde\ref{sphericalweights}.

\begin{lemma}\label{lemma_weights}
Let $[v]\in(V/\mathfrak g.v_{\underline\lambda})^{G_{v_{\underline\lambda}}}$
be of weight $\gamma$.
\smallbreak
\noindent
$(1)$\enspace
If $\gamma=2\alpha$ with $\alpha\in S$,
there exists a unique $j$ such that $(\lambda_j,\alpha)\neq0$.
Moreover, $(\lambda_j,\alpha^\vee)$ is even.
\smallbreak
\noindent
$(2)$\enspace
If $\gamma=\alpha+\beta$ with $\alpha,\beta\in S$ and $(\alpha,\beta)=0$,
there exists a unique $j$ such that
$(\lambda_j,\alpha).(\lambda_j,\beta)\neq0$.
Moreover, $(\lambda_j,\alpha^\vee)=(\lambda_j,\beta^\vee)$
and $(\lambda_i,\alpha)=(\lambda_i,\beta)=0$ for all $i\neq j$.
\end{lemma}

\begin{proof}
(1) Since $v\in\oplus_i V(\lambda_i)_{\lambda_i-2\alpha}$,
one of the dominant weights $\lambda_i$, say $\lambda_j$,
is such that $(\lambda_j,\alpha^\vee)\geq 2$.
Recall that the weights of
$(V/\mathfrak{g}. v_{\underline\lambda})^{G_{v_{\underline\lambda}}}$
are integral linear combinations of the weights $\lambda_i$.
Since the weight of $[v]$ is equal up to a scalar to a simple root,
the coefficients involved in its writing as linear combination
of the weights $\lambda_i$ are all negative except that of $\lambda_j$.
Because of saturation, we thus get:
$(\lambda_i,\alpha)\neq 0$ if and only if $i=j$.
Moreover, we have that $(\lambda_j,\alpha^\vee)$ divides
$(2\alpha,\alpha^\vee)=4$.
With the above inequality this implies that $(\lambda_j,\alpha^\vee)$ is even.

(2) The given simple roots $\alpha,\beta$ generate a root system of type
$\mathsf A_1\times \mathsf A_1$.
Since $v\in\oplus_i V(\lambda_i)_{\lambda_i-\gamma}$,
there exists $\lambda_j$ such that
$(\lambda_j,\alpha).(\lambda_j,\beta)\neq 0$.
By saturation, $\lambda_j$ is the unique dominant weight
which satisfies this inequality.
Therefore, the coordinates $v_i\in V(\lambda_i)$ of $v$ are all $0$
except $v_j$.
The only two pairs which can fulfill condition (\eqref{pair}) are
$(\alpha,\beta)$ and $(\beta,\alpha)$.
Suppose $0\neq X_\alpha.v=X_{-\beta}.v_{\underline\lambda}$,
we then have: $0=X_\alpha.v_i=X_{-\beta}.v_{\lambda_i}$ for $i\neq j$.
From
$0\neq(\lambda_j,\beta^\vee)v_{\lambda_j}=X_\beta X_\alpha.v_j=
X_\alpha X_\beta.v_j=(\lambda_j,\alpha^\vee)v_{\lambda_j}$,
we get also that $X_\beta.v_j\neq0$,
hence $(\lambda_i,\alpha)=(\lambda_i,\beta)=0$ for all $i\neq j$
and $(\lambda_j,\alpha^\vee)=(\lambda_j,\beta^\vee)$.
\end{proof}

\begin{corollary}
\label{ppties_weights}
The set $\Sigma$ of $T_{\ad}$-weights of
$(V/\mathfrak{g}. v_{\underline\lambda})^{G_{v_{\underline\lambda}}}$
satisfies the following properties.
\smallbreak
\noindent
$(\Sigma1)$\enspace
Let $\alpha$ be a simple root such that $2\alpha\in\Sigma$. Then
$(\gamma,\alpha^\vee) \in 2\mathbb Z_{\leq 0}$
for all $\gamma\in\Sigma\setminus\{2\alpha\}$.
\smallbreak
\noindent
$(\Sigma2)$\enspace
Let $\alpha$ and $\beta$ be orthogonal simple roots.
If $\alpha+\beta\in\Sigma$
then $(\gamma,\alpha^\vee)=(\gamma,\beta^\vee)$ for every $ \gamma\in\Sigma$.
\end{corollary}

\begin{proof}
Suppose there exists a weight $\gamma\in\Sigma\setminus\{2\alpha\}$
such that $(\gamma,\alpha)>0$.
By the above argument, we have: $(\gamma,\alpha^\vee)\geq 2$.
Hence in case $\gamma$ is a root,
so are $\gamma-\alpha$ and $\gamma-2\alpha$,
but according to the list given in Theorem\texttilde\ref{sphericalweights},
this can not be.
Consider now, the case where $\gamma$ is not a root.
Checking again the list of possible $T_{\ad}$-weights,
one obtains as single possibility
$\gamma=\alpha_{i-1}+2\alpha_i+\alpha_{i+1}$ with $\alpha_i=\alpha$.
Hence we have in this case a dominant weight
which is not orthogonal to $\alpha_{i-1}$.
Applying Proposition~\ref{insidesupp} with the simple root $\alpha_{i-1}$,
this yields a contradiction.
The first assertion follows.

The second assertion is contained in Lemma\texttilde\ref{lemma_weights}(2).
\end{proof}

\paragraph{Proof of Theorem~\ref{sphericalsystem}.}
Since by Lemma~\ref{pairlemma} the set $\Sigma$ does not contain any simple 
root we need to prove only axioms ($\Sigma1$), ($\Sigma2$) and (S).
The first two follow from Corollary\texttilde\ref{ppties_weights} and
the latter from Proposition\texttilde\ref{insidesupp}
and Remark\texttilde\ref{combinatorialaxiom}. The axiom (S) for the spherical 
root with support of type $\mathsf F$ follows from Lemma~\ref{f4}.

\subsection{Regularity  of $\Gamma$}

Let us first introduce the the set of \textit{colours} $\Delta$ associated
to a spherical system $(S^p,\Sigma,\mathbf A)$ following~\cite{L01}.

The set $\Delta$ is the disjoint union
\begin{equation}
\Delta=\mathbf A\cup\Delta^{a'}\cup\Delta^b  
\end{equation}
where
$\Delta^{a'}$ (resp. $\Delta^b$) contains an element $D_\alpha$ for each $\alpha\in S$ such that $2\alpha\in\Sigma$ 
(resp. $\alpha\not\in S^p\cup\Sigma\cup\frac{1}{2}\Sigma$).
If $\alpha+\alpha'\in\Sigma$ or
$\frac{1}{2}(\alpha+\alpha')\in\Sigma$ then we put $D_\alpha=D_{\alpha'}$.
Note that $\Delta$ is indeed a disjoint union because of the conditions $(\Sigma_1)$ and $(\Sigma_2)$ of spherical systems.

Denote $\omega_\alpha$ the fundamental weight associated to a simple root $\alpha$.
Then one defines a map $\sigma:\Z\Delta\to\Lambda$ as follows (\emph{see}~\cite{F})
\begin{eqnarray}
\label{eq:picbis}
D\mapsto \left \{
\begin{array} {ll}
\sum \omega_\alpha & \quad \mbox{ if  $D\in\mathbf A\cup\Delta^b$ and
  $D=D_\alpha$} \\ 
2\omega_\alpha & \quad  \mbox{ if $D\in\Delta^{a'}$ and $D=D_\alpha$}
\end{array}
\right.
.
\end{eqnarray}


\begin{proposition}\label{char}
Consider the spherical system associated to a free and saturated monoid
$\Gamma=\langle \lambda_1,\ldots,\lambda_s\rangle$.
Let $\Delta$ be its set of colours.
\smallbreak
{\rm(i)}\enspace
We can regard the cone $\mathbb N\Delta$ generated by
the set of colours as a subset of $\Lambda^+$.
The cone $\Gamma$ is then included in the cone $\mathbb N\Delta$.
\smallbreak
{\rm(ii)}\enspace
Each colour of $\Delta$ occurs in the linear combination of
at least one weight $\lambda_i$.
\end{proposition}

The second assertion of the proposition can be called \emph{property of
  regularity}.  

\begin{proof}
Since the set $\mathbf A$ of the spherical system is empty,
the map $\sigma$ defined above is injective.

The image of a colour in $\Lambda^+$ is 
equal to
$2\omega_\alpha$ for all simple roots $\alpha$ with $2\alpha\in\Sigma$,
$\omega_\alpha+\omega_\beta$ for all orthogonal simple roots $\alpha,\beta$
with $\alpha+\beta\in\Sigma$ and
$\omega_\alpha$ for the remaining simple roots $\alpha\notin S^p$.
Since $S^p=\{\alpha\in S:(\lambda_i,\alpha)=0 \,\forall i\}$
the assertion follows from Lemma\texttilde\ref{lemma_weights}.
\end{proof}

\section{The affine multicone over a wonderful variety}\label{sct:multicone}

\subsection{The spherical system of a wonderful variety}

Let $X$ be a wonderful $G$-variety. The set $\Sigma_X$ of its spherical
roots is a basis of the lattice $\Xi_X$ of weights of $B$-semi-invariant
rational functions on $X$. 

A colour of $X$ is a $B$-stable not $G$-stable prime divisor in $X$, namely
the closure of a colour (a $B$-stable prime divisor) of the open $G$-orbit. Let
$\Delta_X$ denote the set of colours of $X$.

If $P_X$ denotes tha standard parabolic subgroup
equal to the stabiliser of the colours of $X$, the closed orbit of $X$ is
isomorphic to $G/P_X$. Let $S^p_X$ denote the set of simple roots
associated to $P_X$.

Let $\mathbf A_X$ denote the set of colours $D$ such that there is a simple
root $\alpha$ in $\Sigma_X$ with $D$ not $P_\alpha$-stable.

Define the map $\rho_X\colon\Delta_X\to\Xi_X^\ast$ by
$\langle\rho_X(D),\gamma\rangle=v_D(f_\gamma)$, where $v_D$ is the
valuation associated to the divisor $D$ and $f_\gamma$ is the
$B$-semi-invariant rational function of weight $\gamma$ (uniquely
determined up to a scalar). 

Then the triple $S^p_X$, $\Sigma_X$, $\mathbf A_X$ is a spherical system
for $G$ and $\Delta_X$ is the set of colours associated to it (\cite{L01}
Section~7). 

\subsection{The open orbit}

Recall the notation of Paragraph~\ref{par:sphericalsystem}.
Let $(S^p,\Sigma,\emptyset)$ be the spherical system for $G$
attached to $\Gamma$.

We shall assume in the following that there exists an associated wonderful
subgroup $N$ of $G$,
that is the generic stabiliser
of a wonderful $G$-variety $X$ with given spherical system $(S^p,\Sigma,\emptyset)$.

In Section\texttilde\ref{sct:existence} we will prove the existence
of such a wonderful subgroup for all spherical systems
arising from a free and saturated cone.

The main purpose of this section is to prove the second assertion of
Theorem~\ref{main} under this existence assumption, that is

\begin{proposition}\label{correspondence}
There exists a subgroup $H$ of $N$, spherical in $G$, such that
the homogeneous space $G/H$ is quasi-affine with weight monoid $\Gamma$. 
We may regard $G/H$ in $V$ then its closure $X_1$ in $V$ is a deformation
of $X_0$.  
\end{proposition}

We need a reformulation of Proposition~\ref{char} in terms of regular
functions on $G$.
Therefore, we shall recall an equivalent definition of the map $\sigma$
given in (\ref{eq:picbis}). 
We consider here the opposite Borel subgroup $B^-$ of $B$. We can suppose
that the set $B^-\,N$ is open in $G$.

The inverse image of a $B^-$-stable prime divisor of $G/N$ through the
quotient map $G\to G/N$ admits 
an equation in $\C[G]$ which is $B^-$-left-semi-invariant and
$N$-right-semi-invariant or for short a $B^-\times N$-eigenvector of $\mathbb
C[G]$.  
This defines an injective map on the set of colours of $G/N$
\begin{equation}\label{eq:pic}
\sigma'\colon\Z\Delta\hookrightarrow\Lambda\times\Xi(N)
\end{equation}
where $\Xi(N)$ is the lattice of $N$-weights.

Consider now the composition of the map $\sigma'$ of (\ref{eq:pic}) with
the projection onto $\Lambda$. 
Such a map
\begin{equation}
\sigma'':\Z\Delta\to\Lambda
\end{equation}
is indeed the map $-\sigma$ (\emph{see} \cite{F}).

We thus get the reformulation of Proposition~\ref{char}.

\begin{proposition}\label{reformulation}
\smallbreak
{\rm(i)}\enspace
For any weight $\chi$ in $\Gamma$, there exists an unique $B^-\times
N$-eigenvector of $\mathbb C[G]$ of $B^-$-weight $-\chi$.
\smallbreak
{\rm(ii)}\enspace
Let $f_i$ be the $B^-\times N$-eigenvector of $\mathbb C[G]$ of weight $-\lambda_i$ for $i=1,\ldots,s$.
Then the function $\Pi_{j=1,\ldots,s} f_j$
vanishes on $G\setminus B^-\,N$.
\end{proposition}

\paragraph{Proof of Proposition~\ref{correspondence}}
Via the isomorphism of $G$-modules
$\C[G]\cong\oplus_{\lambda\in\Lambda^+}V(\lambda)^\ast\otimes V(\lambda)$
the function $f_i$ corresponds to
$\eta_i\otimes v_i\in V(\lambda_i)^\ast\otimes V(\lambda_i)$,
where $\eta_i$ is a $B^-$-eigenvector 
and $v_i$ a $N$-eigenvector.  

Consider the point $x=([v_1],\ldots,[v_s])$ in the product of the projective spaces 
$\P(V(\lambda_i))$. It is therefore $N$-stable.
Further, from the second assertion of
Proposition\texttilde\ref{reformulation}, we have: $N\subset G_x
\subset N_G(N)$. 

We can suppose $N=N_G(N)$.
Indeed, the normaliser $N_G(N)$ acts on the set of colours
of $G/N$ preserving the associated functionals on the lattice $\Xi$
generated by $\Sigma$.
Since in our case there are no spherical roots equal to simple roots,
there can not be distinct colours with the same functional,
hence the spherical homogeneous space $G/N_G(N)$
has the same set of colours as $G/N$
and, by Proposition~\ref{R},
its wonderful compactification has the same spherical system
as that of $G/N$ (\textit{see} \cite{L01}).

Then we have: $G_x=N$ and
\begin{equation}
G/N\ \subset\ \P(V(\lambda_1))\times\ldots\times\P(V(\lambda_s)).
\end{equation}
Define $H$ to be the stabiliser $G_v$ of $v=v_1+\ldots+v_s$ in $V$.
We obtain a spherical homogeneous space, namely $G/H$, with
\begin{equation}
G/H\ \subset\ V.
\end{equation}


Recall that a $N$-invariant $B$-eigenvector of $\mathbb C (G)$,
is uniquely determined by its $B$-weight, up to a scalar. 
Moreover, it is (since $\C[G]$ is a unique factorisation domain) a Laurent
monomial of irreducible $B\times H$-eigenvectors of $\mathbb C [G]$. 
Namely, each spherical root can be seen as an integral combination of colours
via the map $\sigma'$ defined above. 

Therefore, the $B$-weights of the $H$-invariant $B$-eigenvectors
of $\mathbb C [G]$ are in correspondence with the given monoid
$\Gamma$, since it is saturated and $\mathbb Z\Gamma$ contains $\Sigma$.

The closure $X_1$ of $G/H$ in $V$ is thus a deformation of $X_0$ 
(\emph{see} Paragraph~\ref{par:onsat}).


\subsection{The wonderful variety}

Once we know the existence of the above wonderful subgroup $N$,
to obtain Theorem~\ref{main} is enough to show that the closure of 
$G/N$ in the product of the projective spaces $\P(V(\lambda_i))$
is the wonderful embedding. We shall make use of the following result.

\begin{theorem}[\cite{Pe05}]\label{thm:pe}
Let $G$ be (connected) semisimple and simply connected,
let $\overline G$ be its adjoint group $G/Z(G)$.
A wonderful $\overline G$-variety $X$ admits a (unique) closed immersion
in $\P(V(\lambda))$ if
\smallbreak\noindent$(i)$ the dominant weight $\lambda$ is
linear combination with nonzero positive integral coefficients
of the colours of $X$ via the map $\sigma$ in (\ref{eq:picbis}),
\smallbreak\noindent$(ii)$ for any spherical root $\gamma$ of $X$,
there exists no wonderful $\overline G$-variety $X'$ of rank one
with $2\gamma$ as spherical root and $S^p_{X'}=S^p_X$.
\end{theorem}

Take $\lambda$ to be equal to $\lambda_1+\ldots+\lambda_s$.
The closed immersion of $G/N$ in
$\P(V(\lambda_1+\ldots+\lambda_s))$
factorises trough the map
\begin{equation}
\phi\colon\P(V(\lambda_1))\times\ldots\times\P(V(\lambda_s))
\to\P(V(\lambda_1+\ldots+\lambda_s)),
\end{equation}
the Segre embedding followed by the map induced by the projection
of $\otimes_1^s V(\lambda_i)$ onto $V(\lambda_1+\ldots+\lambda_s)$.
The closure of $G/N$ in $\prod_1^s\P(V(\lambda_i))$ dominates
the closure of $G/N$ in $\P(V(\lambda_1+\ldots+\lambda_s))$.
Therefore, if the latter is wonderful then they are equal
for the universal property of wonderful varieties.

Since the first condition of the previous theorem is already fulfilled by
Proposition~\ref{char}, 
to conclude that the closure of $G/N$ in $\prod_1^s\P(V(\lambda_i))$
is wonderful we are left to prove the following.

\begin{proposition}\label{R}
Let $S^p$ and $\Sigma$ defined as above
from a free and saturated cone $\Gamma\subseteq\Lambda^+$.
Then for any $\gamma\in\Sigma$
there exists no wonderful $\overline G$-variety $X'$ of rank one
with $2\gamma$ as spherical root and $S^p_{X'}=S^p$.
\end{proposition}

\begin{proof}
The list of possible elements of $\Sigma$ is in Table~\ref{tab:weights}
and the list of wonderful varieties of rank one is in \cite{W}.
The only case of $\gamma\in\Sigma$ such that $2\gamma$ is a spherical root
is of support of type $\mathsf B_n$ for $\gamma=\alpha_1+\ldots+\alpha_n$
(where $\alpha_1,\ldots,\alpha_n$ are simple roots,
namely the simple roots of the support of $\gamma$).
Let $X'$ be a wonderful variety with $2\alpha_1+\ldots+2\alpha_n$
as spherical root,
then $\alpha_n\in S^p_{X'}$ from Proposition~\ref{insidesupp}.
Suppose $\alpha_1+\ldots+\alpha_n\in\Sigma$ and $\alpha_n\in S^p$,
then there exists a unique $j$ such that $(\lambda_j,\gamma)\neq0$,
it is exactly equal to $\omega_{\alpha_1}$
(the fundamental weight associated to $\alpha_1$).
Hence, if we restrict the action of the torus on
$(V(\lambda_j)/\mathfrak g.v_{\lambda_j})^{G_{v_{\lambda_j}}}$
we should obtain $\alpha_1+\ldots+\alpha_n$ as weight,
but this is false and the weight is actually $2\alpha_1+\ldots+2\alpha_n$.
It can be computed
by taking the semisimple part of the Levi subgroup $G'$ of $G$
with root subsystem generated by $\alpha_1,\ldots,\alpha_n$,
with $T'=G'\cap G$ as maximal torus,
and considering the $T'$-weight in
$(V(\omega'_1)/\mathfrak g'.v')^{G'_{v'}}$
($\omega'_1$ first fundamental weight of $G'$ with respect to $T'$,
$v'$ highest weight vector in $V(\omega'_1)$).
\end{proof}

\begin{corollary}\label{final}
Let ($S^p$, $\Sigma$, $\emptyset$) be a spherical system satisfying the
condition of Proposition~\ref{R}, i.e.\ for any $\gamma\in\Sigma$ there
exists no wonderful variety $X'$ of rank one with $2\gamma$ as spherical
root and $S^p_{X'}=S^p$. 
Let $\Gamma=\langle\lambda_1,\ldots,\lambda_s\rangle_{\mathbb N}$ be free and
saturated, with $S^p$ equal to the subset of simple roots orthogonal to
$\Gamma$ and $\mathbb Z\Gamma\supset\Sigma$, 
satisfying the conditions of Proposition~\ref{char}, i.e.\
$\Gamma\subseteq\mathbb N\Delta$, via the map $\sigma$, and each colour
$D\in\Delta$ occurs in the linear combination of at least one weight
$\lambda_i$.

Assume moreover that, for a fixed $\Gamma$, the set $\Sigma$ is maximal
with the above conditions.

Then ($S^p$, $\Sigma$, $\emptyset$) is the spherical system attached to
$\Gamma$ in the sense of Theorem~\ref{sphericalsystem} and the dimension of
the invariant Hilbert scheme $\Hilb_\Gamma^G(V)_0$ equals the rank of
$\Sigma$.  
\end{corollary}

Notice that the hypotheses of the corollary are fulfilled by $\Gamma=\mathbb
N\Delta$.

\begin{proof}
In the proof of Theorem~\ref{main} we actually use only the relations
between $\Gamma$ and ($S^p$, $\Sigma$, $\emptyset$) that are here reported
as hypotheses of the statement.
\end{proof}

\section{Existence of the wonderful subgroup}\label{sct:existence}

\subsection{Reduction to the primitive cases}


In the previous section we remarked that the spherical system arising from a free and saturated cone $\Gamma$ satisfies the condition of Proposition~\ref{R}. We should now prove that for every spherical system with such a condition there exists a corresponding wonderful subgroup. We want to reduce the proof to some primitive spherical system as in \cite{L01} and \cite{BP}. To do so we use a weaker requirement on the spherical system: we admit all spherical roots which are not simple.

\begin{proposition}\label{existence} Let $S^p$ be a subset of simple roots
  and let $\Sigma$ be a set of spherical not simple roots such that ($S^p$,
  $\Sigma$, $\emptyset$) is a spherical system for $G$ of adjoint
  type. Then there exists a wonderful subgroup $N$ of $G$ such
  that the wonderful embedding of $G/N$ has the given spherical
  system. 
\end{proposition}


In the following we briefly recall some definitions and statements from \cite{L01} which allow us to reduce the proof to finitely many cases. We restrict to the hypotheses of Proposition~\ref{existence}.

Let $(S^p,\Sigma,\emptyset)$ be a spherical system, let $\Delta$ be its set of colours. A subset $\Delta'\subset\Delta$ is said to be \textit{distinguished} if there exists a linear combination $\phi$ of the elements of $\rho(\Delta')$, with positive coefficients, such that $\langle\phi,\gamma\rangle\geq 0$ for all $\gamma\in\Sigma$. Let $\Sigma(\Delta')$ denote the maximal subset of spherical roots such that there exists a linear combination $\phi$ as above such that $\langle\phi,\gamma\rangle>0$ for all $\gamma\in\Sigma(\Delta')$ and let $V(\Delta')$ denote the subspace spanned by $\rho(\Delta')$ and $\{\gamma^\ast\colon\gamma\in\Sigma(\Delta')\}$, where $\{\gamma^\ast\}$ is the dual basis of $\Sigma$.

Let $(S^p,\Sigma,\emptyset)$ be a spherical system and $\Delta'$ a distinguished subset of $\Delta$. Let $\Xi/\Delta'$ be the subgroup of the weights $\xi\in\Xi$ such that $\langle\rho(D),\xi\rangle=0$ for all $D$ in $\Delta'$ and $\langle\gamma^\ast,\xi\rangle=0$ for all $\gamma$ in $\Sigma(\Delta')$. Let $S^p/\Delta'$ be the set $\{\alpha\in S^p\colon \Delta(\alpha)\subset\Delta'\}\subset S^p.$ Let $\Sigma/\Delta'$ be the set of the indecomposable elements of the semigroup $\{\sum_{\gamma\in\Sigma}c_\gamma \gamma$ $\in\Xi/\Delta'\colon$ $c_\gamma\geq0$ $\forall\gamma\in\Sigma\}$. If the set $\Sigma/\Delta'$ is a basis of $\Xi/\Delta'$, then ($S^p/\Delta'$, $\Sigma/\Delta'$, $\emptyset$) is a spherical system.

Let $\Phi\colon X\to X'$ be a dominant $G$-morphism between wonderful
$G$-varieties and let $\Delta_\Phi$ be the set of colours that map dominantly.

\begin{proposition}[\cite{L01}]\label{morphisms} The map $\Phi\mapsto\Delta_\Phi$ is a bijection between the set of dominant $G$-morphisms with connected fibers of $X$, onto another wonderful $G$-variety, and the set of distinguished subsets of $\Delta_X$ such that ($S^p/\Delta_\Phi$, $\Sigma/\Delta_\Phi$, $\emptyset$) is a spherical system. The latter is equal to the spherical system of the target of $\Phi$.\end{proposition}

The distinguished subset $\Delta'$ of $\Delta$ is said to be \textit{smooth} if $V(\Delta')$ is spanned by $\{\gamma^\ast\colon \gamma\in\Sigma(\Delta')\}$. In this case $\Sigma/\Delta'=\Sigma\setminus\Sigma(\Delta')\subset\Sigma$. The distinguished subset $\Delta'$ of $\Delta$ is said to be \textit{parabolic} if $\Sigma(\Delta')=\Sigma$.

Let $P$ be a parabolic subgroup of $G$ and let $L\subset P$ be a Levi subgroup. Let $Y$ be a wonderful $L$-variety. The \textit{parabolic induction} obtained from $Y$ by $P$ is the wonderful $G$-variety $G\times_P Y$, where $Y$ is considered as a $P$-variety where the radical $P^r$ of $P$ acts trivially.

\begin{proposition}[\cite{L01}]\label{P4} Let $X$ be a wonderful $G$-variety and let $S'$ be a subset of $S$ with $(\supp(\Sigma)\cup S^p)\subset S'$. Then $X$ is obtained by parabolic induction by $P^-_{S'}$, the parabolic subset containing $B^-$ associated to $S'$. The corresponding $G$-morphism $\Phi\colon X\to G/P^-_{S'}$ is associated to the parabolic distinguished subset $\Delta(S')$.
\end{proposition}

A spherical system is said to be \textit{cuspidal} if $\supp(\Sigma)=S$.

A spherical system $(S^p,\Sigma,\emptyset)$ is said to be \textit{reducible} if there exists a partition of $S$ into two subsets $S_1,S_2$ such that $S_1 \perp S_2$ and for all $\gamma$ in $\Sigma$ $\supp(\gamma)\subset S_1$ or $\supp(\gamma)\subset S_2$. In this case $(S^p,\Sigma,\emptyset)$ is the \textit{direct product} of the spherical systems ($S^p\cap S_1$, $\Sigma_1$, $\emptyset$) and $(S^p\cap S_2, \Sigma_2, \emptyset)$, where $\Sigma_i=\{\gamma\in\Sigma\colon \supp(\gamma)\subset S_i\}$.

Let $(S^p, \Sigma, \emptyset)$ be a spherical system. Let $\Delta_1$ and $\Delta_2$ be two distinguished subsets of $\Delta$, obviously $\Delta_3=\Delta_1\cup \Delta_2$ is distinguished. The subsets $\Delta_1$ and $\Delta_2$ \textit{decompose into fiber product} the spherical system $(S^p,\Sigma,\emptyset)$ if: \smallbreak\noindent\rm{(i)}\enspace $\Delta_1\neq\emptyset$, $\Delta_2\neq\emptyset$ and $\Delta_1 \cap \Delta_2 = \emptyset$, \smallbreak\noindent\rm{(ii)}\enspace $(S^p/\Delta_i,\Sigma/\Delta_i,\emptyset)$ is a spherical system, for $i=1,2,3$, \smallbreak\noindent\rm{(iii)}\enspace $(\Sigma\setminus(\Sigma/\Delta_1)) \cap (\Sigma\setminus(\Sigma/\Delta_2))=\emptyset$, \smallbreak\noindent\rm{(iv)}\enspace $((S^p/\Delta_1)\setminus S^p)\perp ((S^p/\Delta_2)\setminus S^p)$, \smallbreak\noindent\rm{(v)}\enspace $\Delta_1$ or $\Delta_2$ is smooth.

\begin{proposition} [\cite{L01}] Let $(S^p,\Sigma,\emptyset)$ be a spherical system, let $\Delta_1$ and $\Delta_2$ be two distinguished subsets that decompose $(S^p,\Sigma,\emptyset)$. Let us suppose that, for $i=1,2,3$, there exists a wonderful $G$-variety $X_i$ with spherical system $(S^p/\Delta_i,\Sigma/\Delta_i,\emptyset)$. Then there exists two $G$-morphisms $\Phi_1\colon X_1\to X_3$, $\Phi\colon X_2\to X_3$ such that the fiber product $X_1\times_{X_3}X_2$ is a wonderful $G$-variety with spherical system $(S^p,\Sigma,\emptyset)$.\end{proposition}

The above operation of factorization is a particular case of this decomposition into fiber product. A reducible spherical system, $S=S_1\sqcup S_2$, is decomposed by $\Delta(S_1)$ and $\Delta(S_2)$, and it corresponds to a direct product of two wonderful varieties.

\subsection{Primitive cases}

With the above results the proof of Proposition~\ref{existence} is reduced to the spherical systems that are cuspidal and indecomposable. They are called \textit{primitive} spherical systems.

To list the primitive spherical systems we proceed as follows (\textit{see} \cite{BP}).

Let $(S^p,\Sigma,\emptyset)$ be a spherical system and let $\Delta$ be its set of colours. Two spherical roots $\gamma_1,\gamma_2\in\Sigma$ are said to be \textit{strongly $\Delta$-adjacent} if for all $D\in\Delta(\supp\gamma_1)$ we have $\langle\rho(D),\gamma_2\rangle\neq0$ and, vice versa, for all $D\in\Delta(\supp\gamma_2)$ we have $\langle\rho(D),\gamma_1\rangle\neq0$.

Let $\Sigma^\prime$ be a subset of $\Sigma$. The spherical system $(S^p\cap\supp(\Sigma^\prime),\Sigma^\prime,\emptyset)$ is said to be \textit{strongly $\Delta$-connected} if for every couple of spherical roots $\gamma_1,\gamma_2\in\Sigma^\prime$ there exists a finite sequence of spherical roots in $\Sigma'$, one strongly $\Delta$-adjacent to the next one, the first equal to $\gamma_1$ and the last equal to $\gamma_2$. If $\Sigma^\prime\subset\Sigma$ is maximal with this property we say that $(S^p\cap\supp(\Sigma^\prime),\Sigma^\prime,\emptyset)$ is a \textit{strongly $\Delta$-connected component} of $(S^p,\Sigma,\emptyset)$.

Let $(S^p,\Sigma,\emptyset)$ be a spherical system. Let $\Sigma^\prime$ be a subset of spherical roots. Let $\Delta(\Sigma^\prime)$ denote the subset of colours $D\in\Delta(\supp(\Sigma^\prime))$ such that $\langle\rho(D),\gamma\rangle=0$ for all $\gamma\in\Sigma\setminus\Sigma^\prime$. We say that  $(S^p\cap\supp(\Sigma^\prime),\Sigma^\prime,\emptyset)$ is \textit{erasable} if there exists a nonempty smooth distinguished subset of colours $\Delta^\prime$ included in $\Delta(\Sigma^\prime)$.

We say that it is \textit{quasi-erasable} if there exists a nonempty distinguished subset of colours $\Delta^\prime$ included in $\Delta(\Sigma^\prime)$ such that ($S^p/\Delta'$, $\Sigma/\Delta'$, $\emptyset$) is a spherical system.

\begin{lemma}[\cite{BP}] Let $(S^p,\Sigma,\emptyset)$ be a spherical system. Let $\Sigma_1$ and $\Sigma_2$ be two disjoint subsets of $\Sigma$ giving two quasi-erasable localizations such that at least one of them is erasable. Then the corresponding subsets $\Delta^\prime_1\subset\Delta(\Sigma_1)$ and $\Delta^\prime_2\subset\Delta(\Sigma_2)$ decompose the spherical system $(S^p,\Sigma,\emptyset)$.\end{lemma}

The strongly $\Delta$-connected component $(S^p\cap\supp(\Sigma^\prime),\Sigma^\prime,\emptyset)$ is said to be \textit{isolated} if the partition $\supp(\Sigma^\prime)\sqcup (\supp(\Sigma)\setminus \supp(\Sigma^\prime))$ gives a factorization of $(S^p\cap\supp(\Sigma),\Sigma,\emptyset)$. An isolated component is erasable.


To obtain the list of primitive spherical systems we start from the list of cuspidal strongly $\Delta$-connected spherical systems. Looking at the set of colours of a strongly $\Delta$-connected component we can say whether it is necessarily isolated, erasable, quasi-erasable or none of them. To get all the remaining primitive spherical systems we glue together, in all possible ways, two or more strongly $\Delta$-connected components.


To conclude the proof of Proposition~\ref{existence} we provide here the
list of primitive spherical systems with the corresponding wonderful
subgroups $N$. Since the lists for rank one and two are already
contained in \cite{W}, we restrict to rank greater than two. The proof that
$N$ is wonderful and corresponds to the given spherical system is
already given in \cite{L01}, \cite{Pe03} and \cite{BP} for
types $\mathsf A$, $\mathsf C$, $\mathsf D$. For the other
cases the argument is similar.

The unipotent radical $N^u$ of $N$ is equal to the unipotent radical
$(N^\circ)^u$ of the connected component $N^\circ$. A Levi factor $K$ of
$N$ is such that $N^\circ=K^\circ N^u$. There exists a parabolic subgroup $P$ with Levi decomposition $L\,P^u$ with $K^\circ\subseteq L$ and $N^u\subseteq P^u$. We can suppose that $P$ is in standard position. To describe $N$ we provide here the subgroup $K$ up to isogeny and, if $N$ is not reductive, the subset $S'$ of simple roots associated to $P$ and the complementary $\mathfrak m$ of the Lie algebra of $N^u$ in the Lie algebra $\mathfrak n$ of $P^u$ as representation of the semisimple part $(K^\circ)'$ of $K^\circ$. Moreover, if necessary, we specify the submodule of $\mathfrak n/[\mathfrak n,\mathfrak n]$ in which $\mathfrak m$ is embedded diagonally.

In the following the rank of the connected components of the Dynkin diagram is denoted by $n_i$, or $n$ if there is only one component.

\subsubsection*{Type $\mathsf A$}

\primitivecase $n\geq7$ odd, $S^p=\{\alpha_1,\ldots,\alpha_{2i-1},\alpha_{2i+1},\dots,\alpha_n\}$, $\Sigma=$ $\{\alpha_1+2\alpha_2+\alpha_3,\ldots$, $\alpha_{2i-1}+2\alpha_{2i}+\alpha_{2i+1},\ldots$, $\alpha_{n-2}+2\alpha_{n-1}+\alpha_n\}$. $K=N(Sp_{n+1})$.
\primitivecase $n_1=n_2\geq3$, $S^p=\emptyset$, $\Sigma=\{\alpha_1+\alpha'_1,\ldots$, $\alpha_i+\alpha'_i,\ldots$, $\alpha_{n_1}+\alpha'_{n_2}\}$. $K=N(SL_{n_1+1})$.
\primitivecase $n=2p+q$, $p\geq2$, $q\geq2$, $S^p=\{\alpha_{p+2},\ldots,\alpha_{p+q-1}\}$, $\Sigma=\{\alpha_1+\alpha_{p+q+1},\ldots,\alpha_i+\alpha_{p+q+i},\ldots,\alpha_p+\alpha_n;\alpha_{p+1}+\ldots+\alpha_{p+q}\}$. $K=S(GL_{p+q}\times GL_{p+1})$.
\primitivecase $n=2p+1$, $p\geq2$, $S^p=\emptyset$, $\Sigma=\{\alpha_1+\alpha_{p+2},\ldots,\alpha_i+\alpha_{p+i+1},\ldots,\alpha_p+\alpha_n;2\alpha_{p+1}\}$. $K=N(SL_{p+1}\times SL_{p+1})$.
\primitivecase $n\geq3$, $S^p=\emptyset$, $\Sigma=\{2\alpha_1,\ldots,2\alpha_i,\ldots,2\alpha_n\}$. $K=N(SO_{n+1})$.
\primitivecase $n\geq4$, $S^p=\emptyset$, $\Sigma=\{\alpha_1+\alpha_2,\ldots,\alpha_i+\alpha_{i+1},\ldots,\alpha_{n-1}+\alpha_n\}$. $n$ even, $S'=S\setminus\{\alpha_n\}$, $K=Sp_n\times GL_1$. $n$ odd, $S'=S\setminus\{\alpha_1,\alpha_n\}$, $K=Sp_{n-1}\times GL_1$, $\mathfrak m=V(\omega_1)\subset V(\omega_1)^{\oplus 2}$.

\subsubsection*{Type $\mathsf B$}

\primitivecase $n_1=n_2\geq3$, $S^p=\emptyset$, $\Sigma=\{\alpha_1+\alpha'_1,\ldots$, $\alpha_i+\alpha'_i,\ldots$, $\alpha_{n_1}+\alpha'_{n_2}\}$. $K=N(SO_{2n_1+1})$.
\primitivecase $n\geq4$, $S^p=\emptyset$, $\Sigma=\{\alpha_1+\alpha_2,\ldots,\alpha_i+\alpha_{i+1},\ldots,\alpha_{n-1}+\alpha_n\}$. $n$ even, $S'=S\setminus\{\alpha_n\}$, $K=Sp_n\times GL_1$, $\mathfrak m=V(\omega_1)$. $n$ odd, $S'=S\setminus\{\alpha_1,\alpha_n\}$, $K=Sp_{n-1}\times GL_1$, $\mathfrak m=V(\omega_1)\subset V(\omega_1)^{\oplus 2}$.
\primitivecase $n=p+q$, $p\geq2$, $q\geq1$, $S^p=\{\alpha_{p+2},\ldots,\alpha_n\}$, $\Sigma=\{2\alpha_1,\ldots,2\alpha_i,\ldots$, $2\alpha_p$; $2\alpha_{p+1}+\ldots+2\alpha_n\}$. $K=SO_{p+1}\times SO_{2n-p}$.
\primitivecase $n_1=n_2=2$, $S^p=\emptyset$, $\Sigma=\{\alpha_1+\alpha_2,\alpha_2+\alpha'_2,\alpha'_1+\alpha'_2\}$. $K=SL_2\times SL_2\times SL_2$.
\primitivecase $n=p+q$, $p\geq3$, $q\geq1$, $S^p=\alpha_{p+2},\ldots,\alpha_n$, $\Sigma=\{\alpha_1+\alpha_2,\ldots,\alpha_i+\alpha_{i+1},\ldots,\alpha_{p-1}+\alpha_p;2\alpha_{p+1}+\ldots+2\alpha_n\}$. $p$ even, $S'=S\setminus\{\alpha_p\}$, $K=Sp_p\times N_{SO_{2n-2p+1}}(SO_{2n-2p})\times GL_1$, $\mathfrak m=V(\omega_1)$. $p$ odd, $S'=S\setminus\{\alpha_1,\alpha_p\}$, $K=Sp_{p-1}\times N_{SO_{2n-2p+1}}(SO_{2n-2p})\times GL_1$, $\mathfrak m=V(\omega_1)\subset V(\omega_1)^{\oplus2}$.
\primitivecase $n\geq3$, $S^p=\emptyset$, $\Sigma=\{\alpha_1+\alpha_2,\ldots,\alpha_i+\alpha_{i+1},\ldots,\alpha_{n-1}+\alpha_n;2\alpha_n\}$. $K=N(GL_n)$.
\primitivecase $n=p+q$, $p\geq3$, $q\geq2$, $S^p=\{\alpha_{p+2},\ldots,\alpha_n\}$, $\Sigma=\{\alpha_1+\alpha_2,\ldots,\alpha_i+\alpha_{i+1},\ldots,\alpha_{p-1}+\alpha_p;\alpha_{p+1}+\ldots+\alpha_n\}$. $p$ even, $S'=S\setminus\{\alpha_p\}$, $K=Sp_p\times SO_{2n-2p}\times GL_1$, $\mathfrak m=V(\omega_1)$. $p$ odd, $S'=S\setminus\{\alpha_1,\alpha_p\}$, $K=Sp_{p-1}\times SO_{2n-2p}\times GL_1$, $\mathfrak m=V(\omega_1)\subset V(\omega_1)^{\oplus2}$.

\subsubsection*{Type $\mathsf C$}

\primitivecase $n\geq6$ even, $S^p=\{\alpha_1,\ldots,\alpha_{2i-1},\alpha_{2i+1},\dots,\alpha_{n-1}\}$, $\Sigma=$ $\{\alpha_1+2\alpha_2+\alpha_3,\ldots$, $\alpha_{2i-1}+2\alpha_{2i}+\alpha_{2i+1},\ldots$, $\alpha_{n-3}+2\alpha_{n-2}+\alpha_{n-1};2\alpha_{n-1}+2\alpha_n\}$. $K=N(Sp_n\times Sp_n)$.
\primitivecase $n_1=n_2\geq3$, $S^p=\emptyset$, $\Sigma=\{\alpha_1+\alpha'_1,\ldots$, $\alpha_i+\alpha'_i,\ldots$, $\alpha_{n_1}+\alpha'_{n_2}\}$. $K=N(Sp_2n)$.
\primitivecase $n\geq3$, $S^p=\emptyset$, $\Sigma=\{2\alpha_1,\ldots,2\alpha_i,\ldots,2\alpha_n\}$. $K=N(GL_n)$.
\primitivecase $n=p+q$, $p\geq4$, $q\geq2$, $S^p=\{\alpha_1,\ldots$, $\alpha_{2i-1},\alpha_{2i+1},\dots$, $\alpha_{p+1}$; $\alpha_{p+3},\ldots$, $\alpha_n\}$, $\Sigma=$ $\{\alpha_1+2\alpha_2+\alpha_3,\ldots$, $\alpha_{2i-1}+2\alpha_{2i}+\alpha_{2i+1},\ldots$, $\alpha_{p-1}+2\alpha_p+\alpha_{p+1}$; $\alpha_{p+1}+2\alpha_{p+2}+\ldots+2\alpha_{n-1}+\alpha_n\}$. $K=Sp_{2p+2}\times Sp_{2n-2p-2}$.
\primitivecase $n=p+q$, $p\geq4$, $q\geq1$, $S^p=\{\alpha_3,\ldots,\alpha_{p-2};\alpha_{p+2},\ldots,\alpha_n\}$, $\Sigma=\{\alpha_1+\alpha_p,\alpha_2+\ldots+\alpha_{p-1},\alpha_{p+1}+2\alpha_{p+2}+\ldots+2\alpha_{n-1}+\alpha_n\}$. $S'=S\setminus\{\alpha_{p-1}\}$, $K=SL_p\times Sp_{2n-2p-2}\times GL_1$, $\mathfrak m=V(\omega'_1)$.
\primitivecase $n_1,n_2\geq3$, $S^p=\{\alpha_3,\ldots,\alpha_{n_1};\alpha'_3,\ldots,\alpha'_{n_2}\}$, $\Sigma=\{\alpha_1+\alpha'_1,\alpha_1+2\alpha_2+\ldots+2\alpha_{n_1-1}+\alpha_{n_1},\alpha'_1+2\alpha'_2+\ldots+2\alpha'_{n_2-1}+\alpha'_{n_2}\}$. $K=SL_2\times Sp_{2n_1-2}\times Sp_{2n_2-2}$.
\primitivecase $n=p+q$, $p\geq3$, $q\geq1$, $S^p=\{\alpha_{p+2},\ldots,\alpha_n\}$, $\Sigma=\{\alpha_1+\alpha_2,\ldots,\alpha_i+\alpha_{i+1},\ldots,\alpha_{p-1}+\alpha_p;\alpha_p+2\alpha_{p+1}+\ldots+2\alpha_{n-1}+\alpha_n\}$. $p$ even, $S'=S\setminus\{\alpha_1\}$, $K=Sp_p\times Sp_{2n-p-2}\times GL_1$, $\mathfrak m=V(\omega_1)$. $p$ odd, $S'=S\setminus\{\alpha_1\}$, $K=Sp_{p-1}\times Sp_{2n-p-1}\times GL_1$, $\mathfrak m=V(\omega'_1)$.

\subsubsection*{Type $\mathsf D$}

\primitivecase $n\geq6$ even, $S^p=\{\alpha_1,\ldots,\alpha_{2i-1},\alpha_{2i+1},\dots,\alpha_{n-1}\}$, $\Sigma=$ $\{\alpha_1+2\alpha_2+\alpha_3,\ldots$, $\alpha_{2i-1}+2\alpha_{2i}+\alpha_{2i+1},\ldots$, $\alpha_{n-3}+2\alpha_{n-2}+\alpha_{n-1};2\alpha_n\}$. $K=N(GL_n)$.
\primitivecase $n=p+q$, $p\geq2$, $q\geq2$, $S^p=\{\alpha_{p+2},\ldots,\alpha_n\}$, $\Sigma=\{2\alpha_1,\ldots,2\alpha_i,\ldots$, $2\alpha_p$; $2\alpha_{p+1}+\ldots+2\alpha_{n-2}+\alpha_{n-1}+\alpha_n\}$. $K=S(O_{p+1}\times O_{2n-p-1})$.
\primitivecase $n\geq7$ odd, $S^p=\{\alpha_1,\ldots,\alpha_{2i-1},\alpha_{2i+1},\dots,\alpha_{n-2}\}$, $\Sigma=$ $\{\alpha_1+2\alpha_2+\alpha_3,\ldots$, $\alpha_{2i-1}+2\alpha_{2i}+\alpha_{2i+1},\ldots$, $\alpha_{n-4}+2\alpha_{n-3}+\alpha_{n-2};\alpha_{n-1}+\alpha_{n-2}+\alpha_n\}$. $K=GL_n$.
\primitivecase $n_1=n_2\geq4$, $S^p=\emptyset$, $\Sigma=\{\alpha_1+\alpha'_1,\ldots$, $\alpha_i+\alpha'_i,\ldots$, $\alpha_{n_1}+\alpha'_{n_2}\}$. $K=N(SO_{2n_1})$.
\primitivecase $n\geq4$, $S^p=\emptyset$, $\Sigma=\{2\alpha_1,\ldots,2\alpha_i,\ldots$, $2\alpha_n\}$. $K=N(SO_n\times SO_n)$.
\primitivecase $n=4$, $S^p=\emptyset$, $\Sigma=\{\alpha_1+\alpha_2+\alpha_3,\alpha_3+\alpha_2+\alpha_4,\alpha_4+\alpha_2+\alpha_1\}$. $K=N(G_2)$.
\primitivecase $n\geq4$, $S^p=\emptyset$, $\Sigma=\{\alpha_1+\alpha_2,\ldots,\alpha_i+\alpha_{i+1},\ldots,\alpha_{n-2}+\alpha_n;\alpha_{n-2}+\alpha_n\}$. $n$ even, $S'=S\setminus\{\alpha_1,\alpha_{n-1},\alpha_n\}$, $K=Sp_{n-2}\times GL_1$, $\mathfrak m=V(\omega_1)\subset V(\omega_1)^{\oplus3}$. $n$ odd, $S'=S\setminus\{\alpha_{n-1},\alpha_n\}$, $K=Sp_{n-1}\times GL_1$, $\mathfrak m=V(\omega_1)\subset V(\omega_1)^{\oplus2}$.
\primitivecase $n=p+q$, $p\geq3$, $q\geq2$, $S^p=\{\alpha_{p+2},\ldots,\alpha_n\}$, $\Sigma=\{\alpha_1+\alpha_2,\ldots,\alpha_i+\alpha_{i+1},\ldots,\alpha_{n-2}+\alpha_n;2\alpha_{p+1}+\ldots+2\alpha_{n-2}+\alpha_{n-1}+\alpha_n\}$. $p$ odd, $S'=S\setminus\{\alpha_1,\alpha_p\}$, $K=Sp_{p-1}\times SO_{2q-1}\times GL_1$, $\mathfrak m=V(\omega_1)\subset V(\omega_1)^{\oplus2}$. $p$ even, $S'=S\setminus\{\alpha_p\}$, $K=Sp_p\times SO_{2q-1}\times GL_1$, $\mathfrak m=V(\omega_1)$.

\subsubsection*{Type $\mathsf E$}

\primitivecase $n=7$ (or $8$), $S^p=\{\alpha_2,\alpha_3,\alpha_4,\alpha_5\}$, $\Sigma=\{2\alpha_1+2\alpha_3+2\alpha_4+\alpha_2+\alpha_5,2\alpha_6+2\alpha_5+2\alpha_4+\alpha_2+\alpha_3,2\alpha_7,(2\alpha_8)\}$. $K=N(E_6\times GL_1),E_7\times SL_2$.
\primitivecase $n=7$, $S^p=\{\alpha_2,\alpha_5,\alpha_7\}$, $\Sigma=\{2\alpha_1,2\alpha_3,\alpha_2+2\alpha_4+\alpha_5,\alpha_5+2\alpha_6+\alpha_7\}$. $K=Spin_{12}\times SL_2$.
\primitivecase $n_1=n_2=6,7,8$, $S^p=\emptyset$, $\Sigma=\{\alpha_1+\alpha'_1,\ldots$, $\alpha_i+\alpha'_i,\ldots$, $\alpha_{n_1}+\alpha'_{n_2}\}$. $K=N(E_{n_1})$.
\primitivecase $n=6$, $S^p=\emptyset$, $\Sigma=\{\alpha_1+\alpha_6,\alpha_3+\alpha_5,2\alpha_2,2\alpha_4\}$. $K=SL_6\times SL_2$.
\primitivecase $n=6,7,8$, $S^p=\emptyset$, $\Sigma=\{2\alpha_1,\ldots,2\alpha_i,\ldots,2\alpha_n\}$. $K=Sp_8,SL_8,Spin_{16}$.
\primitivecase $n=6,7,8$, $S^p=\emptyset$,
$\Sigma=\{\alpha_1+\alpha_3,\alpha_2+\alpha_4,\alpha_3+\alpha_4,\ldots$,
$\alpha_i+\alpha_{i+1},\ldots$, $\alpha_{n-1}+\alpha_n\}$. $n=6$, $S'=S\setminus\{\alpha_2\}$, $K=Sp_6\times GL_1$, $\mathfrak m=V(\omega_1)$. $n=7$, $S'=S\setminus\{\alpha_2,\alpha_7\}$, $K=Sp_6\times GL_1$, $\mathfrak m=V(\omega_1)\subset V(\omega_1)^{\oplus2}$. $n=8$, $S'=S\setminus\{\alpha_2\}$, $K=Sp_8\times GL_1$, $\mathfrak m=V(\omega_1)$.
\primitivecase $n=8$, $S^p=\{\alpha_2,\alpha_3,\alpha_4,\alpha_5\}$, $\Sigma=\{2\alpha_1+2\alpha_3+2\alpha_4+\alpha_2+\alpha_5,2\alpha_6+2\alpha_5+2\alpha_4+\alpha_2+\alpha_3,\alpha_7+\alpha_8\}$. $S'=S\setminus\{\alpha_7\}$, $K=F_4\times SL_2\times GL_1$, $\mathfrak m=V(\omega'_1)$.
\primitivecase $n=6$, $S^p=\emptyset$, $\Sigma=\{\alpha_1+\alpha_6,\alpha_3+\alpha_5,\alpha_2+\alpha_4\}$. $S'=S\setminus\{\alpha_4\}$, $K=N_{SL_3\times SL_3}(SL_3)\times SL_2\times GL_1$, $\mathfrak m=V(\omega'_1)$.
\primitivecase $n=7$, $S^p=\{\alpha_2,\alpha_5,\alpha_7\}$,
$\Sigma=\{\alpha_1+\alpha_3,\alpha_2+2\alpha_4+\alpha_5,\alpha_5+2\alpha_6+\alpha_7\}$.
$S'=S\setminus\{\alpha_3\}$, $K=SL_2\times Sp_6\times GL_1$, $\mathfrak m=V(\omega_1)$.

\subsubsection*{Type $\mathsf F$}

\primitivecase $n_1=n_2=4$, $S^p=\emptyset$, $\Sigma=\{\alpha_1+\alpha'_1,\alpha_2+\alpha'_2,\alpha_3+\alpha'_3,\alpha_4+\alpha'_4\}$. $K=F_4$.
\primitivecase $n=4$, $S^p=\emptyset$, $\Sigma=\{2\alpha_1,2\alpha_2,2\alpha_3,2\alpha_4\}$. $K=Sp_6\times SL_2$.
\primitivecase $n=4$, $S^p=\emptyset$, $\Sigma=\{\alpha_1+\alpha_2,2\alpha_3,2\alpha_4\}$. $S'=S\setminus\{\alpha_2\}$, $K=SL_2\times SO_3\times GL_1$, $\mathfrak m=V(\omega_1)$.

\bigbreak\noindent In type $\mathsf G$ there are no primitive spherical systems with rank greater than two.


\section*{}
\begin{small}
Dipartimento di Matematica, Universit\`a di Padova, via G.B.\ Belzoni 7,
35131 Padova, Italy\\
\verb|bravi@math.unipd.it|\\
\\
Mathematisches Institut, Universit\"at zu K\"oln, Weyertal Str.\ 86-90, 50931
K\"oln, Germany\\
\verb|scupit@mi.uni-koeln.de|
\end{small}
\end{document}